\documentclass[11pt]{amsart}
\usepackage{amssymb,mathrsfs}
\title{Langlands Functoriality Conjecture}

\begin{document}

\author{Jae-Hyun Yang}
%\date{July 7, 2007}
%{\textit{To the memory of my mother}}

\address{Department of Mathematics, Inha University,
Incheon 402-751, Korea}
%\newline
%Present address\,: Department of Mathematics, Inha University,
%Incheon 402-751, Korea}
\email{jhyang@inha.ac.kr }

%\begin{document}

\newtheorem{theorem}{Theorem}[section]
\newtheorem{lemma}{Lemma}[section]
\newtheorem{proposition}{Proposition}[section]
\newtheorem{remark}{Remark}[section]
\newtheorem{definition}{Definition}[section]

\renewcommand{\theequation}{\thesection.\arabic{equation}}
\renewcommand{\thetheorem}{\thesection.\arabic{theorem}}
\renewcommand{\thelemma}{\thesection.\arabic{lemma}}
\newcommand{\BR}{\mathbb R}
\newcommand{\BQ}{\mathbb Q}
\newcommand{\BT}{\mathbb T}
\newcommand{\BM}{\mathbb M}
\newcommand{\bn}{\bf n}
\def\charf {\mbox{{\text 1}\kern-.24em {\text l}}}
\newcommand{\BC}{\mathbb C}
\newcommand{\BZ}{\mathbb Z}
\newcommand{\MK}{{\mathcal M}_k}

\thanks{\noindent{Subject Classification:} 11F70, 11R39, 22E50, 22E55.\\
\indent Keywords and phrases: automorphic representations,
automorphic $L$-functions, local Langlands conjectures, global
Langlands conjecture, Langlands functoriality conjecture, converse
theorems}
%\\ \indent This work was supported by 00000}

%\maketitle

\begin{abstract}
{Functoriality conjecture is one of the central and influential
subjects of the present day mathematics. Functoriality is the
profound lifting problem formulated by Robert Langlands in the
late 1960s in order to establish nonabelian class field theory. In
this expository article, I describe the Langlands-Shahidi method,
the local and global Langlands conjectures and the converse
theorems which are powerful tools for the establishment of
functoriality of some important cases, and survey the interesting
results related to functoriality conjecture. }
\end{abstract}

\maketitle

%{\bf Keywords}: Siegel modular variety, fundamental domain, Siegel modular forms, toroidal compactification
%\tableofcontents
\newcommand\tr{\triangleright}
\newcommand\al{\alpha}
\newcommand\be{\beta}
\newcommand\g{\gamma}
\newcommand\gh{\Cal G^J}
\newcommand\G{\Gamma}
\newcommand\de{\delta}
\newcommand\e{\epsilon}
\newcommand\z{\zeta}
\newcommand\vth{\vartheta}
\newcommand\vp{\varphi}
\newcommand\om{\omega}
\newcommand\p{\pi}
\newcommand\la{\lambda}
\newcommand\lb{\lbrace}
\newcommand\lk{\lbrack}
\newcommand\rb{\rbrace}
\newcommand\rk{\rbrack}
\newcommand\s{\sigma}
\newcommand\w{\wedge}
\newcommand\fgj{{\frak g}^J}
\newcommand\lrt{\longrightarrow}
\newcommand\lmt{\longmapsto}
\newcommand\lmk{(\lambda,\mu,\kappa)}
\newcommand\Om{\Omega}
\newcommand\ka{\kappa}
\newcommand\ba{\backslash}
\newcommand\ph{\phi}
\newcommand\M{{\Cal M}}
\newcommand\bA{\bold A}
\newcommand\bH{\bold H}
\newcommand\D{\Delta}

\newcommand\Hom{\text{Hom}}
\newcommand\cP{\Cal P}

\newcommand\cH{\Cal H}

\newcommand\pa{\partial}

\newcommand\pis{\pi i \sigma}
\newcommand\sd{\,\,{\vartriangleright}\kern -1.0ex{<}\,}
\newcommand\wt{\widetilde}
\newcommand\fg{\frak g}
\newcommand\fk{\frak k}
\newcommand\fp{\frak p}
\newcommand\fs{\frak s}
\newcommand\fh{\frak h}
\newcommand\Cal{\mathcal}

\newcommand\fn{{\frak n}}
\newcommand\fa{{\frak a}}
\newcommand\fm{{\frak m}}
\newcommand\fq{{\frak q}}
\newcommand\CP{{\mathcal P}_g}
\newcommand\Hgh{{\mathbb H}_g \times {\mathbb C}^{(h,g)}}
\newcommand\BD{\mathbb D}
\newcommand\BH{\mathbb H}
\newcommand\CCF{{\mathcal F}_g}
\newcommand\CM{{\mathcal M}}
\newcommand\Ggh{\Gamma_{g,h}}
\newcommand\Chg{{\mathbb C}^{(h,g)}}
\newcommand\Yd{{{\partial}\over {\partial Y}}}
\newcommand\Vd{{{\partial}\over {\partial V}}}

\newcommand\Ys{Y^{\ast}}
\newcommand\Vs{V^{\ast}}
\newcommand\LO{L_{\Omega}}
\newcommand\fac{{\frak a}_{\mathbb C}^{\ast}}

\vskip 0.51cm
\begin{center}
\textit{\large To the memory of my parents}
\end{center}
\vskip 2cm

\centerline{\Large \bf Table of Contents}

\vskip 1cm $ \qquad\qquad\qquad\qquad\textsf{\Large \ 1.
Introduction}$\vskip 0.01cm
%\par

$\qquad\qquad\qquad\qquad \textsf{\Large\ 2. Automorphic
$L$-Functions }$\vskip 0.01cm
%\par

$ \qquad\qquad\qquad\qquad  \textsf{\Large\ 3. Local Langlands
Conjecture }$\vskip 0.01cm
%\par

$ \qquad\qquad\qquad\qquad \textsf{\Large\ 4. Global Langlands
Conjecture }$\vskip 0.01cm
%\par

$ \qquad\qquad\qquad\qquad  \textsf{\Large\ 5. Langlands
Functoriality}$\vskip 0.01cm
%\par

$ \qquad\qquad\qquad\qquad \textsf{\Large\ Appendix\,: Converse
Theorems }$\vskip 0.01cm
%\par

%$ \qquad\qquad\qquad\qquad  \textsf{\large\ 7. Examples }$\vskip
%0.01cm
%\par

%$ \qquad\qquad\qquad\qquad  \textsf{\large 12. Motives and Siegel
%Modular Forms }$\vskip 0.01cm
%\par

%$ \qquad\qquad\qquad\qquad\textsf{\large 13. Remark on Cohomology of a Shimura Variety }$
%\vskip 0.01cm
%\par

$ \qquad\qquad\qquad\qquad\textsf{\Large\ References }$
%\par

\newpage

%%%%%%%%%%%%%%%%%%%%%%%%%%%%%%%%%%%%%%%%%%%%%%%%%%%%%%%%%%%%%%%%%%%
%
%                             Sec 1  Introduction
%
%%%%%%%%%%%%%%%%%%%%%%%%%%%%%%%%%%%%%%%%%%%%%%%%%%%%%%%%%%%%%%%%%%%
\begin{section}{{\bf Introduction}}
\setcounter{equation}{0} \vskip 0.3cm Functoriality Conjecture or
the Principle of Functoriality is the profound question that was
raised and formulated by Robert P. Langlands in the late 1960s to
establish nonabelian class field theory and its reciprocity law.
Functoriality conjecture describes deep relationships among
automorphic representations on different groups. This conjecture
can be described in a rough form as follows: To every
$L$-homomorphism $\varphi: {}^LH\lrt {}^LG$ between the $L$-groups
of $H$ and $G$ that are quasi-split reductive groups, there exists
a {\it natural} lifting or transfer of automorphic representations
of $H$ to those of $G$. In 1978, Gelbart and Jacquet \cite{GeJ}
established an example of the functoriality for the symmetric
square ${\rm Sym^2}$ of $GL(2)$ using the converse theorem on
$GL(3)$. In 2002 after 24 years, Kim and Shahidi \cite{KiS2}
established the functoriality for the symmetric cube ${\rm Sym}^3$
of $GL(2)$ and thereafter Kim \cite{Ki1} proved the validity of
the functoriality for the symmetric fourth ${\rm Sym}^4$ of
$GL(2)$. These results have led to breakthroughs toward certain
important conjectures in number theory, those of Ramanujan,
Selberg and Sato-Tate conjectures. We refer to \cite{Sha7} for
more detail on applications to the progress made toward to the
conjectures just mentioned. Recently the Sato-Tate conjecture for
an important class of cases related to elliptic curves has been
verified by Clozel, Harris, Shepherd-Barron and Taylor \cite{CHT,
HST, Maz, Tay}. Past ten years the functoriality for the tensor
product $GL(2)\times GL(2)\lrt GL(4)$ by Ramakrishnan \cite{Ram2},
for the tensor product $GL(2)\times GL(3)\lrt GL(6)$ by Kim and
Shahidi \cite{KiS2}, for the exterior square of $GL(4)$ by Kim
\cite{Ki1} and the weak functoriality to $GL(N)$ for generic
cuspidal representations of split classical groups
$SO(2n+1),\,Sp(2n)$ and $SO(2n)$ by Cogdell, Kim,
Piatetski-Shapiro and Shahidi \cite{CKPS1, CKPS2} were established
by applying appropriate converse theorems of Cogdell and
Piatetski-Shapiro \cite{CP1, CP2} to analytical properties of
certain automorphic $L$-functions arising from the
Langlands-Shahidi method. In fact, the Langlands functoriality was
established only for very special $L$-homomorphisms between the
$L$-groups. It is natural to ask how to find a larger class of
certain $L$-homomorphisms for which the functoriality is valid. It
is still very difficult to answer this question.

\vskip 0.012cm The Arthur-Selberg trace formula has also provided
some instances of Langlands functoriality (see \cite{Art2, AC,
Lan8, Lan9}). In a certain sense, it seems that the trace formula
is a useful and powerful tool to tackle the functoriality
conjecture. Nevertheless the incredible power of Langlands
functoriality seems beyond present technology and knowledge. We
refer the reader to \cite{Lan13} for Langlands' comments on the
limitations of the trace formula. Special cases of functoriality
arises naturally from the conjectural theory of endoscopy
(cf.\,\cite{KoS}), in which a comparison of trace formulas would
be used to characterize the internal structure of automorphic
representations of a given group. I shall not deal with the trace
formula, the base change and the theory of endoscopy in this
article. Nowadays local and global Langlands conjectures are
believed to be encompassed in the functoriality (cf.\,\cite{Lan12,
Lan13}). Quite recently Khare, Larsen and Savin \cite{KLS1, KLS2}
made a use of the functorial lifting from $SO(2n+1)$ to $GL(2n)$,
from $Sp(2n)$ to $GL(2n+1)$ and the theta lifting of the
exceptional group $G_2$ to $Sp(6)$ to prove that certain finite
simple groups $PSp_n({\mathbb F}_{\ell^k}),\,G_2({\mathbb
F}_{\ell^k})$ and $SO_{2n+1}({\mathbb F}_{\ell^k})$ with some mild
restrictions appear as Galois groups over $\BQ$.

\vskip 0.012cm This paper is organized as follows. In Section 2,
we review the notion of automorphic $L$-functions and survey the
Langlands-Shahidi method briefly following closely the article of
Shahidi \cite{Sha3}. I would like to recommend to the reader two
lecture notes which were very nicely written by Cogdell \cite{Co}
and Kim \cite{Ki2} for more information on automorphic
$L$-functions. In Section 3, we review the Weil-Deligne group
briefly and formulate the local Langlands conjecture. We describe
the recent results about the local Langlands conjecture for
$GL(n)$ and $SO(2n+1)$. In Section 4, we discuss the global
Langlands conjecture which is still not well formulated in the
number field case. The work on the global Langlands conjecture for
$GL(2)$ over a function field done by Drinfeld was extended by
Lafforgue ten years ago to give a proof of the global Langlands
conjecture for $GL(n)$ over a function field. We will not deal
with the function field case in this article. We refer to
\cite{Laf} for more detail. Unfortunately there is very little
known of the global Langlands conjecture in the number field case.
I have an audacity to mention the Langlands hypothetical group and
the hypothetical motivic Galois group following the line of
Arthur's argument in \cite{Art1}. In Section 5, I formulate the
Langlands functoriality conjecture in several ways and describe
the striking examples of Langlands functoriality established past
ten years. I want to mention that there is a descent method of
studying the opposite direction of the lift initiated by Ginzburg,
Rallis and Soudry (see \cite{GRS1, JiS2}). I shall not deal with
the descent method here. In the appendix, I describe a brief
history of the converse theorems obtained by Hamburger, Hecke,
Weil, Cogdell, Piatetski-Shapiro, Jinag and Soudry. I present the
more or less exact formulations of the converse theorems. As
mentioned earlier, the converse theorems play a crucial role in
establishing the functoriality for the examples discussed in
Section 5.

\vskip 0.1cm \noindent {\bf Notations:} \ \ We denote by
$\BQ,\,\BR$ and $\BC$ the field of rational numbers, the field of
real numbers and the field of complex numbers respectively. We
denote by $\BR^*_+$ the multiplicative group of positive real
numbers. $\BC^*$ (resp. $\BR^*$) denotes the multiplicative group
of nonzero complex (resp. real) numbers. We denote by $\BZ$ and
$\BZ^+$ the ring of integers and the set of all positive integers
respectively. For a number field $F$, we denote by ${\mathbb A}_F$
and ${\mathbb A}_F^*$ the ring of adeles of $F$ and the
multiplicative group of ideles of $F$ respectively. If there is no
confusion, we write simply $\mathbb A$ and ${\mathbb A}^*$ instead
of ${\mathbb A}_F$ and ${\mathbb A}_F^*$. For a field $k$ we
denote by $\Gamma_k$ the Galois group $Gal({\overline k}/k)$,
where $\overline k$ is a separable algebraic closure of $k$. We
denote by ${\mathbb G}_m$ the multiplicative group in one
variable. ${\mathbb G}_a$ denotes the additive group in one
variable.

\end{section}

\vskip 1cm
%%%%%%%%%%%%%%%%%%%%%%%%%%%%%%%%%%%%%%%%%%%%%%%%%%%%%%%%%%%%%%%%%%%
%
%   Sec 14 The Weil Representation
%
%%%%%%%%%%%%%%%%%%%%%%%%%%%%%%%%%%%%%%%%%%%%%%%%%%%%%%%%%%%%%%%%%%%
%\begin{section}{\bf The Weil Representation}
%\setcounter{equation}{0}
\newcommand\A{\mathbb A}
\newcommand\GA{G({\mathbb A})}
\newcommand\PA{P({\mathbb A})}
\newcommand\MA{M({\mathbb A})}
\newcommand\NA{N({\mathbb A})}
\newcommand\KA{K_{\mathbb A}}
\newcommand\FO{\mathfrak o}
\newcommand\BA{B({\mathbb A})}
\newcommand\LG{{}^L G}
\newcommand\LM{{}^L M}
\newcommand\LN{{}^L N}
\newcommand\LLN{{}^L {\mathfrak n}}
\newcommand\LP{{}^L P}

\begin{section}{{\bf  Automorphic $L$-functions}}
\setcounter{equation}{0} \vskip 0.3cm Let $G$ be a connected,
quasi-split reductive group over a number field $F$. For each
place $v$ of $F$, we let $F_v$ be the completion of $F$,
${\mathfrak o}_v$ the rings of integers of $F_v$, ${\mathfrak
p}_v$ the maximal ideal of ${\mathfrak o}_v$, and let $q_v$ be the
order of the residue field $k_v={\mathfrak o}_v /{\mathfrak p}_v$.
We denote by ${\mathbb A}={\mathbb A}_F$ the ring of adeles of
$F$. We fix a Borel subgroup $B$ of $G$ over $F$. Write $B=TU$,
where $T$ is a maximal torus and $U$ is the unipotent radical,
both over $F$. Let $P$ be a parabolic subgroup of $G$. Assume
$P\supset B.$ Let $P=MN$ be a Levi decomposition of $P$ with Levi
factor $M$ and its unipotent radical $N$. Then $N\subset U.$
\vskip 0.13cm For each place $v$ of $F$, we let $G_v=G(F_v).$
Similarly we use $B_v,T_v,U_v,P_v,M_v$ and $N_v$ to denote the
corresponding groups of $F_v$-rational points. Let
$\GA,\BA,\cdots,\NA$ be the corresponding adelic groups for the
subgroups defined before. When $G$ is unramified over a place $v$
in the sense that $G$ is quasi-split over $F_v$ and that $G$ is
split over a finite unramified extension of $F_v$, we let
$K_v=G(\FO_v).$ Otherwise we shall fix a special maximal compact
subgroup $K_v\subset G_v.$ We set $\KA=\otimes_v K_v.$ Then
$\GA=\PA\KA.$

\vskip 0.1cm Let $\Pi=\otimes_v \Pi_v$ be a cuspidal automorphic
representation of $\GA$. We refer to \cite{JaS1, JaS2, Lan6} for
the notion of automorphic representations. Let $S$ be a finite set
of places including all archimedean ones such that both $\Pi_v$
and $G_v$ are unramified for any place $v\notin S$. Then for each
$v\notin S,\ \Pi_v$ determines uniquely a semi-simple conjugacy
class $c(\Pi_v)$ in the $L$-group ${}^LG_v$ of $G_v$ as a group
defined over $F_v$. We refer to \cite{B, KZ1} for the definition
and construction of the $L$-group. We note that there exists a
natural homomorphism $\xi_v:{}^LG_v\lrt {}^LG.$ For a finite
dimensional representation $r$ of $\LG$, putting $r_v=r\circ
\xi_v$, the local Langlands $L$-function $L(s,\Pi_v,r_v)$
associated to $\Pi_v$ and $r_v$ is defined to be
(cf.\,\cite{B},\,\cite{Lan1})
\begin{equation}
L(s,\Pi_v,r_v)=\det \Big( I-r_v\big(c(\Pi_v)\big)
q_v^{-s}\Big)^{-1}.
\end{equation}
We set
\begin{equation}
L_S(s,\Pi,r)=\prod_{v\notin S}L(s,\Pi_v,r_v).
\end{equation}

Langlands (cf.\,\cite{Lan1}) proved that $L_S(s,\Pi,r)$ converges
absolutely for sufficiently large ${\rm Re} (s)> 0$ and defines a
holomorphic function there. Furthermore he proposed the following
question. \vskip 0.1cm \noindent {\bf Conjecture A
(Langlands\,\cite{Lan1}).} $L_S(s,\Pi,r)$ has a meromorphic
continuation to the whole complex plane and satisfies a standard
functional equation. \vskip 0.21cm

F. Shahidi (cf.\,\cite{Sha3},\,\cite{Sha4}) gave a partial answer
to the above conjecture using the so-called Langlands-Shahidi
method. I shall describe Shahidi's results briefly following his
article \cite{Sha3}.

\vskip 0.1cm Let $P$ be a maximal parabolic subgroup of $G$ and
$P=MN$ its Levi decomposition with its Levi factor $M$ and its
unipotent radical $N$. Let $A$ be the split torus in the center of
$M$. For every group $H$ defined over $F$, we let $X(H)_F$ be the
group of $F$-rational characters of $H$. We set
\begin{equation*}
{\mathfrak a}=\Hom \big( X(M)_F,\BR\big).
\end{equation*}
Then
\begin{equation*}
{\mathfrak a}^*=X(M)_F \otimes_\BZ \BR\cong X(A)_F\otimes_\BZ \BR.
\end{equation*}
We set ${\mathfrak a}^*_\BC:={\mathfrak a}^*\otimes_\BR \BC.$ Let
${\mathfrak z}$ be the real Lie algebra of the split torus in the
center $C(G)$ of $G$. Then ${\mathfrak z}\subset {\mathfrak a}$
and ${\mathfrak a}/{\mathfrak z}$ is of dimension $1$. The
imbedding $X(M)_F\hookrightarrow X(M)_{F_v}$ induces an imbedding
$ {\mathfrak a}_v \hookrightarrow {\mathfrak a},$ where
${\mathfrak a}_v=\Hom \big( X(M)_{F_v},\BR\big).$ The
Harish-Chandra homomorphism $H_P:\MA\lrt {\mathfrak a}$ is defined
by
\begin{equation*}
\exp \langle \chi,H_P(m)\rangle=\prod_{v} |\chi(m_v)|_v,\quad
\chi\in X(M)_F,\ m=\otimes_v m_v\in\MA.
\end{equation*}
We may extend $H_P$ to $\GA$ by letting it trivial on $\NA$ and
$\KA$. We define $H_{P_v}:M_v\lrt {\mathfrak a}_v$ by
\begin{equation*}
q_v^{ \langle \chi_v,H_{P_v}(m)\rangle}= |\chi_v(m_v)|_v,\quad
\chi_v\in X(M)_{F_v},\ m_v\in M_v
\end{equation*}
for a finite place $v$, and define
\begin{equation*}
\exp{ \langle \chi_v,H_{P_v}(m)\rangle}= |\chi_v(m_v)|_v,\quad
\chi_v\in X(M)_{F_v},\ m_v\in M_v
\end{equation*}
for an infinite place $v$. Then we have
\begin{equation}
\exp \langle \chi,H_P(m)\rangle=\prod_{v=\infty} \exp{ \langle
\chi,H_{P_v}(m)\rangle}\cdot \prod_{v<\infty} q_v^{ \langle
\chi,H_{P_v}(m)\rangle},
\end{equation}
where $\chi\in X(M)_F$ and $ m=\otimes_v m_v\in\MA.$ We observe
that the product in (2.3) is finite.

\vskip 0.12cm Let $A_0$ be the maximal $F$-split torus in $T$. We
denote by $\Phi$ the set of roots of $A_0.$ Then $\Phi=\Phi^+ \cup
\Phi^-,$ where $\Phi^+$ is the set of roots generating $U$ and
$\Phi^-=-\Phi^+.$ Let $\Delta\subset \Phi^+$ be the set of simple
roots. The unique reduced root of $A$ in $N$ can be identified by
an element $\alpha\in\Delta.$ Let $\rho_P$ be half the sum of
roots generating $N$. We set
\begin{equation*}
{\widetilde \alpha}=\langle \rho_P,\alpha\rangle^{-1} \rho_P.
\end{equation*}
Here, for any pair of roots $\alpha$ and $\beta$ in $\Phi^+$, the
pairing $\langle \alpha,\beta\rangle$ is defined as follows. Let
${\widetilde \Phi}^+$ be the set of non-restricted roots of $T$ in
$U$. We see that the set of simple roots $\widetilde{\Delta}$ in
${\widetilde \Phi}^+$ restricts to $\Delta$. Identifying $\alpha$
and $\beta$ with roots in ${\widetilde \Phi}^+$, we set
$$\langle \alpha,\beta\rangle ={{2(\alpha,\beta)}\over
{(\beta,\beta)}},$$ where $(\,\ ,\ )$ is the standard inner
product on $\BR^l$ with $l=|{\widetilde\Delta}|$.

Let $\pi=\otimes_v \pi_v$ be a cuspidal automorphic representation
of $\MA$. Given a $\KA\cap M(\A)$-finite function $\phi$ in the
representation space $V_\pi$ of $\pi$, we may extend $\phi$ to a
function $\widetilde{\phi}$ on $\GA$\,(cf.\,\cite{Sha1}). Then the
Eisenstein series $E(s,\widetilde{\phi},g,P)$ is defined by
\begin{equation}
E(s,\widetilde{\phi},g,P)=\sum_{\gamma\in P(F)\backslash G(F)}
\widetilde{\phi}(\gamma g)\exp \langle s
\widetilde{\alpha}+\rho_P,H_P(\gamma g)\rangle,\quad g\in \GA.
\end{equation}
The Eisensetein series $E(s,\widetilde{\phi},g,P)$ converges
absolutely for sufficiently large ${\rm Re}(s)\gg 0$ and extends
to a meromorphic function of $s$ on $\BC$, with a finite number of
poles in the plane ${\rm Re}(s)>0$, all simple and on the real
axis (See \cite{Lan5}).

 Let $W$ be the Weyl group of $A_0$ in $G$. We denote the
subset of $\Delta$ which generates $M$ by $\theta$. Then
$\Delta=\theta\cup \{ \alpha\}.$ Then there exists a unique
element $\widetilde{w}\in W$ such that
$\widetilde{w}(\theta)\subset \Delta$ and
$\widetilde{w}(\alpha)\in \Phi^-.$ Fix a representative $w\in
\KA\cap G(F)$ for $\widetilde{w}$. We shall also denote every
component of $w$ by $w$ again.

\vskip 0.12cm Let
\begin{equation*}
I(s,\pi)={\rm Ind}_{\MA\NA}^{\GA} \pi\otimes \exp \langle
s\tilde{\alpha},H_P(\,\cdot\,)\rangle \otimes 1
\end{equation*}
be the representation of $\GA$ induced from $\PA$. Then
$I(s,\pi)=\otimes_v I(s,\pi_v)$ with
\begin{equation*}
I(s,\pi_v)={\rm Ind}_{M_vN_v}^{G_v} \pi_v\otimes q_v^{ \langle
s\tilde{\alpha},H_P(\,\cdot\,)\rangle }\otimes 1,
\end{equation*}
where $q_v$ should be replaced by $\exp$ if $v$ is archimedean. We
let $M'$ be the subgroup of $G$ generated by
$\widetilde{w}(\theta)$. Then there is a parabolic subgroup
$P'\supset B$ with $P'=M'N'$. Here $M'$ is a Levi factor of $P'$
and $N'$ is the unipotent radical of $P'$. For $f\in I(s,\pi)$ and
sufficiently large ${\rm Re} (s)\gg 0,$ we define
\begin{equation}
M(s,\pi)f(g)=\int_{N'(\mathbb A)}f(w^{-1}ng)dn,\quad g\in \GA.
\end{equation}
At each place $v$, for sufficiently large ${\rm Re} (s)\gg 0,$ we
define a local intertwining operator by
\begin{equation}
A(s,\pi_v,w)f_v(g)=\int_{N_v'}f_v(w^{-1}ng)dn,\quad g\in G_v,
\end{equation}
where $f_v\in I(s,\pi_v).$ Then
\begin{equation}
M(s,\pi)=\otimes_v A(s,\pi_v,w).
\end{equation}
It follows from the theory of Eisenstein series that for ${\rm Re}
(s)\gg 0,$ $M(s,\pi)$ extends to a meromorphic function of
$s\in\BC$ with only a finite number of simple poles
(cf.\,\cite{Lan5}).

\vskip 0.12cm Let $\LM$ and $\LN$ be the Levi factor and the
unipotent radical of the parabolic subgroup $\LP=\LM\LN$ of the
$L$-group $\LG.$ Then we have the representation $r:\LM\lrt {\rm
End}\big( \LLN\big)$ given by the adjoint action of $\LM$ on the
Lie algebra $\LLN$ of $\LN.$ Let
$$r=r_1\oplus r_2\oplus\cdots \oplus r_m$$
be the decomposition of $r$ into irreducible constituents. Each
irreducible constituent $(r_i,V_i)$ with $1\leq i\leq m$ is
characterized by
$$V_i=\left\{\, X_{\beta^\vee}\in \LLN\,|\ \langle
\tilde{\alpha},\beta\rangle =i\,\right\},\quad i=1,\cdots,m.$$ We
refer to \cite{Lan2} and \cite[Proposition 4.1]{Sha3} for more
detail.

\vskip 0.12cm According to \cite{Lai} and \cite{Lan2}, one has
\begin{equation}
M(s,\pi)f=\Big( \otimes_{v\in S} A(s,\pi_v,w)f_v\Big) \bigotimes
\left( \otimes_{v\notin S}\tilde{f_v}\right) \times \prod_{i=1}^m
{ {L_S(is,\pi,{r}_i)}\over {L_S(1+is,\pi,{r}_i)} },
\end{equation}
where $f=\otimes_v f_v$ is an element in $I(s,\pi)$ such that for
each $v\notin S,\ f_v$ is the unique $K_v$-fixed vector with
$f_v(e_v)=1,$ $\tilde{f_v}$ is the $K_v$-fixed vector in
$I(-s,\widetilde{w}(\pi_v))$ with $\tilde{f_v}(e_v)=1,$ and
$\widetilde{r}_i$ denotes the contragredient of $r_i\,(1\leq i\leq
m)$.

\vskip 0.12cm For every archimedean place $v$ of $F$, let
$\varphi_v:W_{F_v}\lrt \LM_v$ be the corresponding homomorphism
(cf.\,\cite{Lan10}) attached to $\pi_v$. One has a natural
homomorphism $\eta_v:\LM_v\lrt \LM.$ We put
$$r_{i,v}=r_i\circ \eta_v,\quad i=1,2,\cdots,m.$$
Then $r_{i,v}\circ \varphi_v=r_i\circ \eta_v\circ \varphi_v$ is a
finite dimensional representation of the Weil group $W_{F_v}$ on
$V_i$. Let $L(s,r_{i,v}\circ \varphi_v)$ be the corresponding
Artin $L$-function attached to $r_{i,v}\circ \varphi_v$ (cf.
\,\cite{Lan4}). We set
\begin{equation}
L^S(s,\pi,r_i)=\prod_{v=\infty}L(s,r_{i,v}\circ \varphi_v)\cdot
\prod_{v\notin S\atop v< \infty} L(s,\pi_v,r_i\circ\eta_v).
\end{equation}
Let $\rho:M(\A)\lrt {\overline M}(\A)$ be the projection of
$M(\A)$ onto its adjoint group.

\vskip 0.2cm Shahidi showed the following.

\vskip 0.21cm \noindent {\bf Theorem 2.1\,(Shahidi \cite{Sha3}).}
Let $\pi=\otimes_v\pi_v$ be a cuspidal automorphic representation
of ${\overline M}(\A).$ Then every $L$-function
$L^S(s,\pi,r_i\circ {}^L\rho),\,1\leq i\leq m,$ extends to a
meromorphic function of $s$ to the whole complex plane. Moreover,
if $\pi$ is {\it generic}, then each $L^S(s,\pi,r_i\circ
{}^L\rho)$ satisfies a standard functional equation, that is,
\begin{equation*}
L^S(s,\pi,r_i\circ {}^L\rho)=\varepsilon_S(s,\pi,r_i\circ
{}^L\rho)\, L^S(1-s,\widetilde{\pi,r_i\circ {}^L\rho}),
\end{equation*}
where $\varepsilon_S(1s,\pi,r_i\circ {}^L\rho)$ is the root number
attached to $\pi$ and $r_i\circ {}^L\rho$, and $\widetilde \tau$
denotes the contragredient of a representation $\tau$.

\vskip 0.12cm Furthermore, for a given generic cuspidal
automorphic representation $\pi=\otimes_v \pi_v$ of $M(\mathbb
A)$, Shahidi defined the local $L$-functions $L(s,\pi_v,r_i)$ and
the local root numbers $\varepsilon(s,\pi_v,r_i,\psi)$ with $1\leq
i\leq m$ at {\it bad places} $v$ so that the completed
$L$-function $L(s,\pi,r_i)$ and the completed root number
$\varepsilon(s,\pi,r_i)$ defined by
\begin{equation*}
L(s,\pi,r_i)=\prod_{{\rm all}\ v} L(s,\pi_v,r_i),\quad
\varepsilon(s,\pi,r_i)=\prod_{{\rm all}\ v}
\varepsilon(s,\pi_v,r_i,\psi),\quad i=1,\cdots,m
\end{equation*}
satisfy a standard functional equation
\begin{equation}
L(s,\pi,r_i)= \varepsilon(s,\pi,r_i)\,L(1-s,\pi,{\widetilde
{r_i}}),\quad i=1,\cdots,m.
\end{equation}

\vskip 0.21cm \noindent {\bf Example 2.2\,(Kim-Shahidi
\cite{KiS1}).} Let $F$ be a number field and let $G$ be a simply
connected semisimple split group of type $G_2$ over $F$. We set
${\mathbb A}_{\infty}=\prod_{v=\infty}F_v$. Let $K_{\infty}$ be
the standard maximal compact subgroup of $G({\mathbb A}_{\infty})$
and $K_v=G({\mathfrak o}_v)$ for a finite place $v$. Then
$K_{\mathbb A}=K_{\infty}\times \prod_{v< \infty} K_v$ is a
maximal compact subgroup of $G(\mathbb A).$ Fix a split maximal
torus $T$ in $G$ and let $B=TU$ be a Borel subgroup of $G$. In
what follows the roots are those of $T$ in $U$. Let
$\Delta=\{\beta_1,\beta_6\}$ be a basis of the root system $\Phi$
with respect to $(T,B)$ with the long simple root $\beta_1$ and
the short one $\beta_6$. Then the other roots are given by
$$\beta_2=\beta_1+\beta_6,\quad \beta_3=2\beta_1+3\beta_6,\quad
\beta_4=\beta_1+2\beta_6,\quad \beta_5=\beta_1+3\beta_6.$$ Let $P$
be the maximal parabolic subgroup of $G$ generated by $\beta_1$
with Levi decomposition $P=MN$, where $M\simeq GL(2).$ See
\cite[Lemma 2.1]{Sha6}. Thus one has
\begin{equation*}
{\mathfrak a}^*=\BR \,\beta_4,\quad {\mathfrak a}=\BR\,
\beta_4^{\vee}\quad {\rm and}\quad \rho_P={\frac 52}\,\beta_4.
\end{equation*}
\indent Let ${\widetilde \alpha}=\beta_4.$ Then $s{\widetilde
\alpha}\,(s\in\BC)$ corresponds to the character $|\det(m)|^s$. We
note that ${\mathbb A}^*={\mathbb A}^*_1\cdot \BR^*_+,$ where
${\mathbb A}^*_1$ is the group of ideles of norm $1$. Let
$\pi=\otimes_v\pi_v$ be a cuspidal automorphic representation of
$M({\mathbb A})=GL(2,{\mathbb A}).$ We may and will assume that
the central character $\omega_\pi$ of $\pi$ is trivial on
$\BR_+^*.$ For a $K$-finite function $\varphi$ in the
representation space of $\pi$, the Eisenstein series
$E(s,{\widetilde \varphi},g)=E(s,{\widetilde \varphi},g,P)$
defined by Formula (2.4) converges absolutely for sufficiently
large ${\rm Re}(s)\gg 0$ and extends to a meromorphic function of
$s$ on $\BC$, with a finite number of poles in the plane ${\rm
Re}(s)>0$, all simple and on the real axis. The discrete spectrum
$L^2_{disc}\big(G(F)\backslash G({\mathbb A})\big)_{(M,\pi)}$ is
spanned by the residues of Eisenstein series for ${\rm Re}(s)>0$
(See \cite{Lan5}). We know that the poles of Eisenstein series
coincide with those of its constant terms. So it is enough to
consider term along $P$. For each $f\in I(s,\pi)$, the constant
term of $E(s,f,g)$ along $P$ is given by
$$E_0(s,f,g)=\sum_{w\in\Omega}M(s,\pi,w) f(g),\quad \Omega=\{1,s_6 s_1 s_6 s_1 s_6\},$$
where $s_i$ is the reflection along $\beta_i$ defined by
$$ s_i (\beta)=\beta -{ {2(\beta_i,\beta)}\over
{(\beta_i,\beta_i)}} \beta_i,\quad 1\leq i\leq 6,\ \beta\in
\Phi.$$
Weyl group representatives are all chosen to lie in $K_{\mathbb A}
\cap G(F).$ Here
$$M(s,\pi,w)f(g)=\int_{N_w^-(\mathbb A)} f(w^{-1}ng)dn=\prod_v \int_{N_w^-(F_v)}f_v(w^{-1}_vn_vg_v)dn_v,$$
where $g=\otimes_v g_v\in G(\mathbb A)$, $f=\otimes_v f_v$ is an
element of $I(s,\pi)$ such that $f_v$ is the unique $K_v$-fixed
function normalized by $f_v(e_v)=1$ for almost all $v$, and
$$N_w^-=\prod_{\alpha>0\atop w^{-1}\alpha <0}U_\alpha,\quad
U_\alpha={\rm the\ one\ parameter\ unipotent\ subgroup.}$$ \indent
Let ${\rm St}:GL(2,\BC)\lrt GL(2,\BC)$ be the standard
representation of $GL(2,\BC)$ and
$${\rm Sym}^3:GL(2,\BC)={}^LM\lrt GL(4,\BC)$$
be the third symmetric power representation of $GL(2,\BC).$ Let
$$\big( {\rm Sym}^3\big)^0={\rm Sym}^3 \bigotimes \big(
\wedge^2{\rm St}\big)^{-1}$$ be the adjoint cube representation of
$GL(2,\BC)$ (cf.\,\cite[p.\,249]{Sha6}). Then the adjoint
representation $r$ of ${}^LM=GL(2,\BC)$ on the Lie algebra
${}^L{\mathfrak n}$ of ${}^LN$ is given by
$$r=\big( {\rm Sym}^3\big)^0\oplus \wedge^2 {\rm St}.$$
\indent According to Formula (2.8), one obtain, for $w=s_6 s_1 s_6
s_1 s_6 ,$ the formula
\begin{eqnarray*}
M(s,\pi)f=& &\Big( \otimes_{v\in S} M(s,\pi_v,w)f_v\Big)
\bigotimes
\left( \otimes_{v\notin S}\tilde{f_v}\right) \\
\ & &\times\, L_S\big(s,{\widetilde\pi},( {\rm
Sym}^3)^0\big)\,L_S\big(2s,{\widetilde\pi},\wedge^2 {\rm
St}\big)\\
\ & &\times \, L_S\big(1+s,{\widetilde\pi},( {\rm
Sym}^3)^0\big)^{-1}\,L_S\big(1+2s,{\widetilde\pi},\wedge^2 {\rm
St}\big)^{-1},
%{{L_S\big(s,{\widetilde\pi},( {\rm
%Sym}^3)^0\big)\,L_S\big(2s,{\widetilde\pi},\wedge^2 {\rm St}\big)
%}\over {{L_S\big(1+s,{\widetilde\pi},( {\rm
%Sym}^3)^0\big)\,L_S\big(1+2s,{\widetilde\pi},\wedge^2 {\rm
%St}\big) } }},
\end{eqnarray*}
where $S$ is a finite set of places of $F$ including all the
archimedean places such that $\pi_v$ is unramifies for every
$v\notin S.$ Here $L_S\big(s,\pi,\wedge^2{\rm St}\big)$ is the
partial Hecke $L$-function. Kim and Shahidi \cite{KiS1} proved
that if $\pi$ is a non-monomial cuspidal representation of
$M(\mathbb A)=GL(2,{\mathbb A})$ in the sense that $\pi\not\cong
\pi\otimes\eta$ for any nontrivial grossencharacter $\eta$ of
$F^*\backslash {\mathbb A}_F^+$, the symmetric cube $L$-function
$L\big(s,\pi,{\rm Sym}^3\big)$ and the adjoint cube $L$-function
$L\big(s,\pi,({\rm Sym}^3)^0\big)$ are both {\it entire} and
satisfy the standard functional equations
$$L\big(s,\pi,{\rm Sym}^3\big)=\varepsilon \big(s,\pi,{\rm
Sym}^3\big)\, L\big(1-s,{\widetilde\pi},{\rm Sym}^3\big)$$ and
$$L\big(s,\pi,({\rm Sym}^3)^0\big)=\varepsilon \big(s,\pi,({\rm
Sym}^3)^0\big)\, L\big(1-s,{\widetilde\pi},({\rm Sym}^3)^0\big).$$
It follows from this fact that if $\pi$ is not monomial, the
partial Rankin triple $L$-function $L_S(s,\pi\times\pi\times\pi)$
is entire. Ikeda \cite{Ike} calculated the poles of the Rankin
triple $L$-function $L_S(s,\pi\times\pi\times\pi)$ for $GL(2)$.
And we have the following relations
\begin{equation}
L_S(s,\pi\times\pi\times\pi)=L\big(s,\pi,{\rm
Sym}^3\big)\,\big(L_S(s,\pi\otimes \omega_\pi)\big)^2
\end{equation}
and
\begin{equation}
L\big(s,\pi,{\rm Sym}^3\big)=L_S\big(s,\pi\otimes \omega_\pi,({\rm
Sym}^3)^0\big).
\end{equation}

According to Formula (2.10), $L_S(s,\pi\times\pi\times\pi)$ could
have double zeros at $s=1/2.$ \hfill $\square$

\vskip 0.15cm In \cite{KiS3}, Kim and Shahidi studied the
cuspidality of the symmetric fourth power ${\rm Sym}^4(\pi)$ of a
cuspidal representation $\pi$ of $GL(2,{\mathbb A})$ and the
partial symmetric $m$-th power $L$-functions $L_S\big(s,\pi,{\rm
Sym}^m\big)\,(1\leq m\leq 9).$ For the definition of ${\rm
Sym}^m(\pi)$, we refer to Example 5.6 in this article. We
summarize their results.

\vskip 0.2cm \noindent {\bf Theorem
2.3\,(Kim-Shahidi\,\cite{KiS3}).} Let $\pi$ be a cuspidal
automorphic representation of $GL(2,{\mathbb A})$ with
$\omega_\pi$ its central character. Then ${\rm Sym}^4(\pi)\otimes
\omega_\pi^{-1}$ is a cuspidal representation of $GL(5,{\mathbb
A})$ except for the following three cases:

\vskip 0.15cm \noindent (1) $\pi$ is monomial in the sense that
$\pi\cong \pi\otimes\eta$ for some nontrivial Gr{\"o}ssencharacter
$\eta$ of $F$.\\
(2) $\pi$ is not monomial and ${\rm Sym}^3(\pi)\otimes
\omega_\pi^{-1}$ is not cuspidal.\\
(3) ${\rm Sym}^3(\pi)\otimes \omega_\pi^{-1}$ is cuspidal and
there exists a nontrivial quadratic character $\chi$ such that\\
\indent\  ${\rm Sym}^3(\pi)\otimes \omega_\pi^{-1} \cong {\rm
Sym}^3(\pi)\otimes \omega_\pi^{-1}\otimes \chi.$

\vskip 0.2cm As applications of Theorem 2.3, they obtained the
following.

\vskip 0.2cm \noindent {\bf Proposition
2.4\,(Kim-Shahidi\,\cite{KiS3}).} Let $\pi$ be a cuspidal
automorphic representation of $GL(2,{\mathbb A})$ with
$\omega_\pi$ its central character such that ${\rm Sym}^3(\pi)$ is
cuspidal. Then the following statements hold:

\vskip 0.12cm\noindent (a) Each partial symmetric $m$-th power
$L$-functions $L_S\big(s,\pi,{\rm Sym}^m\big)\,(m=6,7,8,9)$ has
a\\ \indent \ meromorphic continuation and satisfies a standard
functional
equation.\\
(b) $L_S\big(s,\pi,{\rm Sym}^5\big)$ and $L_S\big(s,\pi,{\rm
Sym}^7\big)$ are holomorphic and nonzero for ${\rm Re}(s)\geq
1.$\\
(c) If $\omega_\pi^3=1,$ $L_S\big(s,\pi,{\rm Sym}^6\big)$ is
holomorphic and nonzero for ${\rm Re}(s)\geq
1.$\\
(d) If ${\rm Sym}^4(\pi)$ is cuspidal, $L_S\big(s,\pi,{\rm
Sym}^6\big)$ is holomorphic and nonzero for ${\rm Re}(s)\geq
1.$\\
(e) If ${\rm Sym}^4(\pi)$ is cuspidal and $\omega_\pi^4=1$,
$L_S\big(s,\pi,{\rm Sym}^8\big)$ is holomorphic and nonzero\\
\indent\ for ${\rm Re}(s)\geq
1.$\\
(f) If ${\rm Sym}^4(\pi)$ is cuspidal, $L_S\big(s,\pi,{\rm
Sym}^9\big)$ has a most a simple pole or a simple zero at $s=1.$\\
(g) If ${\rm Sym}^4(\pi)$ is not cuspidal, $L_S\big(s,\pi,{\rm
Sym}^9\big)$ is holomorphic and nonzero for ${\rm Re}(s)\geq 1.$

\vskip 0.2cm \noindent {\bf Proposition
2.5\,(Kim-Shahidi\,\cite{KiS3}).} Let $\pi=\otimes_v\pi_v$ be a
cuspidal automorphic representation of $GL(2,{\mathbb A})$ such
that ${\rm Sym}^3(\pi)$ is cuspidal. Let ${\rm
diag}(\alpha_v,\beta_v)$ be the Satake parameter for $\pi_v$. Then
$|\alpha_v|,|\beta_v| < q_v^{1/9}.$ If ${\rm Sym}^4(\pi)$ is not
cuspidal, then $|\alpha_v|=|\beta_v|=1$, that is, the full
Ramanujan conjecture holds.

\vskip 0.2cm \noindent {\bf Proposition
2.6\,(Kim-Shahidi\,\cite{KiS3}).} Let $\pi=\otimes_v\pi_v$ be a
nonmonomial cuspidal automorphic representation of $GL(2,{\mathbb
A})$ with a trivial central character. Suppose $m\leq 9.$ Then the
following statements hold:

\vskip 0.015cm \noindent (1) Suppose ${\rm Sym}^3(\pi)$ is not
cuspidal. Then $L_S\big(s,\pi,{\rm Sym}^m\big)$ is holomorphic and
nonzero at\\
\indent\
 $s=1$, except for $m=6,8$\,; the $L$-functions
$L_S\big(s,\pi,{\rm Sym}^6\big)$ and $L_S\big(s,\pi,{\rm
Sym}^8\big)$ each\\ \indent \ have a simple pole at $s=1.$\\
(2) Suppose ${\rm Sym}^3(\pi)$ is cuspidal but ${\rm Sym}^4(\pi)$
is not cuspidal. Then $L_S\big(s,\pi,{\rm Sym}^m\big)$ is
holo-\\
\indent \ morphic and nonzero at $s=1$ for $m=1,\cdots,7$ and
$m=9$\,; the $L$-function $L_S\big(s,\pi,{\rm Sym}^8\big)$\\
\indent\ has a simple pole at $s=1.$

\vskip 0.212cm We are still far from solving Conjecture $A$. We
have two known methods to study analytic properties of automorphic
$L$-functions. The first is the method of constructing explicit
zeta integrals that is called the Rankin-Selberg method. The
second is the so-called Langlands-Shahidi method I just described
briefly. In the late 1960s Langlands \cite{Lan2} recognized that
many automorphic $L$-functions occur in the constant terms of the
Eisenstein series associated to cuspidal automorphic
representations of the Levi subgroups of maximal parabolic
subgroups of split reductive groups through his intensively deep
work on the theory of Eisenstein series. He obtained some analytic
properties of certain automorphic $L$-functions using the
meromorphic continuation and the functional equation of Eisenstein
series. As mentioned above, he proved the meromorphic continuation
of certain class of $L$-functions but did not gave an answer to
the functional equation. Shahidi \cite{Sha3} generalized
Langlands' recognition to quasi-split groups, and calculated
non-constant terms of the Eisenstein series and hence obtained the
functional equations of a more broader class of many automorphic
$L$-functions. In fact, those $L$-functions dealt with by Shahidi
include most of the $L$-functions studied by other mathematicians
(cf.\,\cite{GS},\,\cite{GPSR}). The first method has some
advantage to provide more precise information on the location of
poles and the special values of automorphic $L$-functions. On the
other hand, the Langlands-Shahidi method has been applied to a
large class of automorphic $L$-functions, and is likely to be more
suited to the theory of harmonic analysis on a reductive group.
Moreover the second method plays an important role in
investigating the non-vanishing of automorphic $L$-functions on
the line ${\rm Re}(s)=1.$ One of the main contributions of Shahidi
to the Langlands-Shahidi method is to define local $L$-functions
even at bad places in such a way that the completed $L$-function
satisfies the functional equation. We are in need of new methods
to have more knowledge on the analytic and arithmetic properties
of automorphic $L$-functions.

\end{section}

\vskip 1cm
%%%%%%%%%%%%%%%%%%%%%%%%%%%%%%%%%%%%%%%%%%%%%%%%%%%%%%%%%%%%%%%%%%%
%
%   Sec 3  Local Langlands Conjecture
%
%%%%%%%%%%%%%%%%%%%%%%%%%%%%%%%%%%%%%%%%%%%%%%%%%%%%%%%%%%%%%%%%%%%
\begin{section}{{\bf Local Langlands Conjecture}}
\setcounter{equation}{0}

\vskip 0.1cm Let $k$ be a local field and let $W_k$ be its Weil
group. We review the definition of the Weil group $W_k$ following
the article of Tate (cf.\,\cite{Tat}). If $k=\BC$, then
$W_\BC=\BC^\times.$ If $k=\BR,$ then
$$W_\BR=\BC^*\cup \tau \BC^*,\quad \tau z\tau^{-1}=\overline{z},$$
where $z\mapsto \overline{z}$ is the nontrivial element of ${\rm
Gal}(\BC/\BR).$ Then $W_\BR^{\rm ab}=\BR^*.$ Here if $Y^c$ denotes
the closure of the commutator subgroup of a topological group $Y$,
we set $Y^{\rm ab}=Y/Y^c$.

\vskip 0.12cm Suppose $k$ is a nonarchimedean local field and
${\overline k}$ a separable algebraic closure of $k$. Let $q$ be
the order of the residue field $\kappa$ of $k$. We set $\G_k={\rm
Gal}\big({\overline k}/k \big)$ and $\G_{\kappa}={\rm
Gal}\big({\overline \kappa}/\kappa \big)$. Let $\Phi_\kappa:
x\mapsto x^q$ be the Frobenius automorphism in $\G_\kappa.$ We set
$\langle \Phi_\kappa\rangle=\left\{ \Phi_\kappa^n\,|\
n\in\BZ\,\right\}$. Let $\varphi:\G_k\lrt \G_\kappa$ be the
canonical surjective homomorphism given by $\sigma\mapsto
\sigma|_{\overline\kappa}.$ The Weil group $W_k$ is defined to be
the set $W_k=\varphi^{-1}\big( \langle \Phi_\kappa\rangle \big).$
Obviously one has an exact sequence
$$1\lrt I_k\lrt W_k\lrt \langle \Phi_\kappa\rangle\lrt 1,$$
where $I_k={\rm ker}\,\varphi$ is the inertia group of $k$. We
recall that $W_k$ is topologized such that $I_k$ has the induced
topology from $\G_k$, $I_k$ is open in $\G_k$ and multiplication
by $\Phi$ is a homeomorphism. Here $\Phi$ denotes a choice of a
geometric Frobenius element in
$\varphi^{-1}(\Phi_\kappa\big)\subset \G_k.$ We note that we have
a continuous homomorphism $W_k\lrt \G_k$ with dense image.
According to the local class field theory, one has the isomorphism
\begin{equation}
k^*\cong W_k^{\rm ab}.
\end{equation}

\newcommand\Ga{{\mathbb G}_a}
\newcommand\GKG{\mathcal{G}_k(G)}
\newcommand\GKH{\mathcal{G}_k(H)}
\newcommand\PG{\prod(G(k))}
In order to generalize the isomorphism (3.1) for $GL(1)$ to
$GL(2)$, P. Deligne \cite{D2} introduced the so-called {\it
Weil}-{\it Deligne group} $W_k'$. It is defined to be the group
scheme over $\BQ$ which is the semidirect product of $W_k$ by the
additive group $\Ga$ on which $W_k$ acts by the rule
$wxw^{-1}=||w||x.$ We refer to \cite[p.\,19]{Tat} for the
definition of $||w||$. Namely, $W_k'=W_k \ltimes \Ga$ is the group
scheme over $\BQ$ with the multiplication
$$(w_1,a_1)(w_2,a_2)=(w_1w_2, a_1+||w_1||\,a_2),\quad w_1,w_2\in
W_k,\ a_1,a_2\in \Ga.$$

\noindent{\bf Definition 3.1\,(Deligne \cite{D3}).} Let $E$ be
field of characteristic $0$. A representation of $W_k'$ over $E$
is a pair $\rho'=(\rho,N)$ consisting of

\vskip 0.1cm (a) A finite dimensional vector space $V$ over $E$
and a homomorphism $\rho:W_k\lrt GL(V)$ whose kernel contains an
open subgroup of $I_k$, i.e., which is continuous for the discrete
topology on $GL(V)$.\\
\indent (b) A nilpotent endomorphism $N$ of $V$ such that
$$\rho(w)N\rho(w)^{-1}=||w||\,N,\quad w\in W_k.$$

\vskip 0.2cm We see that a homomorphism of group schemes over $E$
$$\rho':W_k\times_\BQ E\lrt GL(V)$$
determines, and is determined by a pair $(\rho,N)$ as in
Definition 3.1 such that
$$\rho'((w,a))=\exp (aN)\,\rho(w),\quad w\in W_k,\ a\in\Ga.$$

Let $\rho'=(\rho,N)$ be a representation of $W_k'$ over $E$.
Define $\nu:W_k\lrt \BZ$ by $||w||=q^{-\nu (w)},\ w\in W_k.$ Then
according to \cite[(8.5)]{D3}, there is a {\it unique\ unipotent}
automorphism $u$ of $V$ such that
$$uN=Nu,\quad u\rho(w)=\rho(w) u,\quad w\in W_k$$
\noindent and such that
$$ \exp (aN)\,\rho(w)\,u^{-\nu(w)}$$

\noindent is a semisimple automorphism of $V$ for all $a\in E$ and
$w\in W_k-I_k.$ Then $\rho'=(\rho u^{-\nu},N)$ is called the
$\Phi$-semisimplication of $\rho'$. And $\rho'$ is called
$\Phi$-semisimple if and only if $\rho'=\rho_{ss}',\ u=1,$ i.e.,
the Frobeniuses acts semisimply.

\vskip 0.15cm Let $\rho'=(\rho,N,V)$ be a representation of $W_k'$
over $E$. We let $V_N^I:=(\ker N)^{I_k}$ be the subspace of
$I_k$-invariants in $\ker N.$ We define a local $L$-factor by
\begin{equation}
Z(t,V)=\det\left( 1-t\,\rho(\Phi)\big|_{V_N^I}\right)^{-1},\ \
{\rm{and}} \quad L(s,V)=Z(q^{-s},V),\ {\rm when}\ E\subset \BC.
\end{equation}

We note that if $\rho'=(\rho,N)$ is a representation of $W_k'$,
then $\rho'$ is irreducible if and only if $N=0$ and $\rho$ is
irreducible.

\vskip 0.2cm Let $G$ be a connected reductive group over a local
field. A homomorphism $\alpha: W_k'\lrt {}^LG$ is said to be {\it
admissible}\,\cite[p.\,40]{B} if the following conditions
(i)-(iii) are satisfied \vskip 0.15cm\noindent
 (i) $\alpha$ is a homomorphism over $\G_k$, i.e., the following
 diagram is commutative:
 \begin{equation*}
\begin{array}{ccc}
W_k' &\stackrel{\alpha}\longrightarrow&  {}^LG\\
\ \ \ \,\,\,\,\, \searrow&&\!\!\!\!\!\!\!\swarrow\\
&  \Gamma_k  &\\
\end{array}
\end{equation*}
 (ii) $\alpha$ is continuous, $\alpha(\Ga)$ are unipotent in
${}^LG^0,$ and $\alpha$ maps semisimple elements into semisimple
elements in ${}^LG$. Here an element $x$ is said to be {\it
semisimple} if $x\in I_k$, and an element $g\in {}^LG$ is called
{\it semisimple} if $r(g)$ is semisimple for any finite
dimensional representation $r$ of ${}^LG$.

(iii) If $\alpha(W_k')$ is contained in a Levi subgroup of a
proper parabolic subgroup $P$ of ${}^LG$, then $P$ is relevant.
See \cite[p.\,32]{B}.

\vskip 0.2cm Let $\GKG$ be the set of all admissible homomorphisms
$\phi:W_k'\lrt \LG$ modulo inner automorphisms by elements of
$\LG^0$. We observe that we can associate canonically to $\phi\in
\GKG$ a character $\chi_\phi$ of the center $C(G)$ of $G$
(cf.\,\cite[p.\,43]{B},\,\cite{Lan10}). Let $Z_L^0=C(\LG^0)$ be
the center of $\LG^0$. Following \cite[pp.\,43-44]{B} and
\cite{Lan10}, we can construct a character $\omega_{\alpha}$ of
$G(k)$ associated to a cohomology class $\alpha\in H^1
(W_k',Z_L^0).$ If we write $\phi\in\GKG$ as $\phi=(\phi_1,\phi_2)$
with $\phi:W_k'\lrt \LG^0$ and $\phi:W_k'\lrt \G_k$, then $\phi_1$
defines a cocycle of $W_k'$ in $\LG^0$, and the map $\phi\mapsto
\phi_1$ yields an embedding $\GKG\hookrightarrow H^1
(W_k',\LG^0).$ Then $H^1 (W_k',Z_L^0)$ acts on $H^1 (W_k',\LG^0)$
and this action leaves $\GKG$ stable \cite[p.\,40]{B}.

\vskip 0.2cm Let $\prod(G(k))$ be the set of all equivalence
classes of irreducible admissible representations of $G(k)$. The
following conjecture gives an arithmetic parametrization of
irreducible admissible representations of $G(k).$

\vskip 0.3cm\noindent {\bf Local Langlands Conjecture [LLC].} Let
$k$ be a local field. Let $\GKG$ and $\prod(G(k))$ be as above.
Then there is a surjective map $\PG \lrt \GKG$ with finite fibres
which partitions $\PG$ into disjoint finite sets
$\prod_\phi(G(k))$, simply $\prod_\phi$ called $L$-{\it packets}
satisfying the following (i)-(v):

\vskip 0.15cm (i) If $\pi\in \prod_\phi$, then the central
character $\chi_\pi$
of $\pi$ is equal to $\chi_\phi.$\\
\indent (ii) If $\alpha\in H^1 (W_k',Z_L^0)$ and $\omega_\alpha$
is its
associated character of $G(k)$, then\\
\begin{center}
%\begin{equation*}
$\prod_{\alpha\cdot\phi}=\left\{ \pi\omega_\alpha\,|\
\pi\in\prod_\phi\ \right\}.$
%\end{equation*}
\end{center}

(iii) One element of $\prod_\phi$ is square integrable modulo the
center $C(G)$ of $G$ if and only if all elements are square
integrable modulo the center $C(G)$ of $G$ if and only if
$\phi(W_k')$ is not contained in any proper Levi subgroup of
$\LG$.\\
\indent (iv) One element of $\prod_\phi$ is tempered if and only
if all elements of $\prod_\phi$ are tempered if and only if
$\phi(W_k)$ is bounded.\\
\indent (v) If $H$ is a connected reductive group over $k$ and
$\eta:H(k)\lrt G(k)$ is a $k$-morphism with commutative kernel and
cokernel, then there is a required compatibilty between
decompositions for $G(k)$ and $H(k)$. More precisely, $\eta$
induces a canonical map ${}^L\eta:\LG\lrt {}^LH,$ and if we set
$\phi'={}^L\eta\circ \phi$ for $\phi\in\GKG,$ then any
$\pi\in\prod_\phi(G(k))$, viewed as an $H(k)$-module, decomposes
into a direct sum of finitely many irreducible admissible
representations belonging to $\prod_{\phi'}(H(k)).$

\vskip 0.3cm \noindent {\bf Remark 3.2.} (a) If $k$ is
archimedean, i.e., $k=\BR$ or $\BC$, [LLC] was solved by Langlands
\cite{Lan10}. We also refer the reader to \cite{ABV, AV, KZ2}.
\\
(b) In case $k$ is non-archimedean, Kazhdan and Lusztig \cite{KaL}
had shown how to parametrize those irreducible admissible
representations of $G(k)$ having an Iwahori fixed vector in terms
of admissible homomorphisms of $W_k'$.\\
(c) For a local field $k$ of positive characteristic $p>0,$ [LLC]
was established by Laumon, Rapoport and Stuhler \cite{NRS}.\\
(d) In case $G=GL(n)$ for a non-archimedean local field $k$, [LLC]
was established by Harris and Taylor \cite{HT}, and by Henniart
\cite{Hen2}. In both cases, the correspondence was established at
the level of a correspondence between irreducible Galois
representations and supercuspidal representations.\\
(e) Let $k$ be a a non-archimedean local field of characteristic
$0$ and let $G=SO(2n+1)$ the split special orthogonal group over
$k$. In this case, Jiang and Soudry \cite{JiS1,JiS2} gave a
parametrization of {\it generic\ supercuspidal} representations of
$SO(2n+1)$ in terms of admissible homomorphisms of $W_k'.$ More
precisely, there is a unique bijection of the set of conjugacy
classes of all admissible, completely reducible,
multiplicity-free, symplectic complex representations
$\phi:W_k'\lrt {}^LSO(2n+1)=Sp(2n,\BC)$ onto the set of all
equivalence classes of irreducible generic supercuspidal
representations of $SO(2n+1,k).$

\vskip 0.3cm For $\pi\in \prod_\phi(G)$ with $\phi\in\GKG,$ if $r$
is a finite dimensional complex representation of $\LG$, we define
the $L$- and $\varepsilon$-factors

\begin{equation}
L(s,\pi,r)=L(s,r\circ \phi)\quad {\rm and}\quad
\varepsilon(s,\pi,r,\psi)=\varepsilon(s,r\circ \phi,\psi),
\end{equation}

\noindent where $L(s,r\circ \phi)$ is the Artin-Weil $L$-function.

\vskip 0.3cm\noindent {\bf Remark 3.3.} For a non-archimedean
local field $k$, Deligne \cite{D2} gave the complete formulation
of [LLC] for $GL(2).$ In \cite{D2}, he utilized for the first time
the Weil-Deligne group $W_k'$, which was introduced by him in
\cite{D3}, in the context of $\ell$-adic representations, in order
to obtain a correct formulation in the case of $GL(2)$ over a
non-archimedean local field.

\vskip 0.3cm\noindent {\bf Remark 3.4.} The representations in the
$L$-packet $\prod_\phi$ are parametrized by the component group
$$C_\phi:=S_\phi/ Z_L S_\phi^0,$$
where $S_\phi$ is the centralizer of the image of $\phi$ in $\LG,\
S_\phi^0$ is the identity component of $S_\phi$, and $Z_L$ is the
center of $\LG$. We refer the reader to \cite{Art3, LL} for more
information on the $L$-packets.

\vskip 0.3cm\noindent {\bf Example 3.5.} Let $\pi$ be a spherical
or unramified representation of $G(k)$ with a non-archimedean
local field $k$. It is known that $\pi\hookrightarrow I(\chi)$ for
a unique unramified quasi-character $\chi$ of a maximal torus
$T(k)$ of $G(k).$ Then $\pi$ determines a semi-simple conjugate
class $c(\pi)=\{ t_\pi\}\subset {}^LT \subset {}^LG.$ Then the
Langlands' parameter $\phi_\pi$ for $\pi$ is

\begin{equation*}
\phi_\pi:k^*\lrt {}^LT\subset {}^LG,\quad\
\phi_\pi(\widetilde{\omega})=t_\pi
\end{equation*}
such that
\begin{equation*}
\phi_\pi(\widetilde{\omega})=t_\pi\quad {\rm and}\ \phi_\pi \ {\rm
is\ trivial\ on}\ {\mathcal O}^*,
\end{equation*}
where $\widetilde{\omega}$ denotes a uniformizer in $k$. \vskip
0.2cm If $\pi=\pi(\mu_1,\cdots,\mu_n)$ is a spherical
representation of $GL_n(k)$ with unramified characters
$\mu_i\,(i=1,\cdots,n)$ of $k^*$, the semi-simple conjugacy class
$c(\pi)$ is given by
$$ c(\pi)=\{ \,{\rm
diag}(\mu_1(\widetilde{\omega}),\cdots,\mu_n(\widetilde{\omega}))\}$$
and the Langlands' parameter $\phi_\pi$ for $\pi$ is
$$ \phi_\pi:k^*\lrt {}^LT\subset {}^LG=GL(n,\BC),\quad\
\phi_\pi(\widetilde{\omega})={\rm
diag}(\mu_1(\widetilde{\omega}),\cdots,\mu_n(\widetilde{\omega})).$$

\end{section}

\vskip 1cm
%%%%%%%%%%%%%%%%%%%%%%%%%%%%%%%%%%%%%%%%%%%%%%%%%%%%%%%%%%%%%%%%%%%
%
%   Sec 4  Global Langlands Conjecture
%
%%%%%%%%%%%%%%%%%%%%%%%%%%%%%%%%%%%%%%%%%%%%%%%%%%%%%%%%%%%%%%%%%%%
\begin{section}{{\bf Global Langlands Conjecture}}
\setcounter{equation}{0}
\newcommand\AuG{{\mathcal A}(G)}
\newcommand\Mk{{\mathcal M}_k}
\newcommand\AKG{{\mathcal A}_k(G)}
\newcommand\LKG{{\mathfrak L}_k(G)}
\newcommand\AKvG{{\mathcal A}_{k_v}(G)}
\newcommand\LKvG{{\mathfrak L}_{k_v}(G)}

\vskip 0.2cm Let $k$ be a global field and ${\mathbb A}$ its ring
of adeles. This section is based on Arthur's article \cite{Art1}.

\vskip 0.012cm As in the local case of Section 3, the global
Langlands conjecture should be a nonabelian generalization of
abelian global class field theory. When Deligne \cite{D3}
recognized the need to introduce the Weil-Deligne group $W_k'$ for
the local Langlands correspondence for $GL(2)$, it was realized
that there seemed to be no natural global version of $W_k'$. In
fact, $\G_k,\,W_k$ and $W_k'$ are too small to parameterize all
automorphic representations of a reductive group. Thus in the
1970s Langlands \cite{Lan7} attempted to discover a {\it
hypothetical} group $L_k$ to replace the Weil-Deligne group
$W_k'$. Nowadays it is believed by experts that this group $L_k$
should be related to the equally {\it hypothetical motivic Galois}
group ${\mathcal M}_k$ of $k$. The group $L_k$ is often called the
{\it hypothetical} (or {\it conjectual}) Langlands group or the
{\it automorphic} Langlands group. The notion of $L_k$ and
${\mathcal M}_k$ is still not clear.

The global Langlands conjecture can be written as follows.

\vskip 0.2cm\noindent {\bf Global Langlands Conjecture [GLC].}
Automorphic representations of $\GA$ can be parametrized by
admissible homomorphisms $\phi:L_k\lrt \LG$ required to have the
following properties (1)-(4): \vskip 0.051cm \noindent  (1) There
is an $L$-packet $\prod_\phi$ which consists
of automorphic representations of $\GA$ attached to $\phi$.\\
\ \ (2) Each $L$-packet $\prod_\phi$ is a finite set.\\
\ \ (3) Any automorphic representation of $\GA$ belongs to
$\prod_\phi$ for a unique homomorphism $\phi.$\\
\ \ (4) The $\prod_\phi$'s are disjoint.

\vskip 0.2cm We first consider the case that $k$ is a function
field. Drinfeld \cite{Dr2} formulated a version of the global
Langlands conjecture for function fields relating the irreducible
two dimensional representations of the Galois group $\G_k$ with
irreducible cuspidal representations of $GL(2,\A)$, and
established the global Langlands conjecture. In the early 2000s
Lafforgue \cite{Laf} had extended the work of Drinfeld mentioned
above to $GL(n)$ to give a proof of the global Langlands
conjecture for $GL(n)$ over a function field. The formulation of
the global Langlands conjecture made by Drinfeld and Lafforgue is
essentially the same as that in the local non-archimedean case
discussed in Section 3 with a few modification. So we omit the
details for the global Langlands conjecture over a function field
here. We refer to \cite{Dr1, Dr2, Laf} for more detail.

\vskip 0.2cm Next we consider the case $k$ is a number field. We
first recall that according to the global class field theory, for
$n=1$, there is a canonical bijection between the continuous
characters of $\G_k$ and characters of finite order of the idele
group $k^*\backslash {\mathbb A}^*$. We should replace $\G_k$ by
the Weil group $W_k$ in order to obtain all the characters of
$k^*\backslash {\mathbb A}^*$. For $n\geq 2,$ by analogy with the
local Langlands conjecture, we need a global analogue of the
Weil-Deligne group $W_k'$. However no such analogue is available
at this moment. We hope that the hypothetical Langlands group
$L_k$ plays a role as $W_k'$ and fits into an exact sequence
\begin{equation}
1 \lrt L_k^c \lrt L_k \lrt \G_k \lrt 1,
\end{equation}
where $L_k^c$ is a complex pro-reductive group. $L_k$ should be a
locally compact group equipped with an embedding $i_v:L_{k_v}\lrt
L_k$ for each completion $k_v$ of $k$. Let $G$ be a connected,
quasi-split reductive group over $k$. We set $G_v:=G(k_v)$ for
every place $v$ of $k$. Let $\LKG$ be the set of all equivalence
classes of continuous, completely reducible homomorphisms $\phi$
of $L_k$ into $\LG,$ and $\AKG$ the set of equivalence classes of
all automorphic representations of $\GA$. For each place $v$ of
$k$, let $\LKvG$ be the set of equivalence classes of continuous,
completely reducible homomorphisms $\phi_v:L_{k_v}\lrt {}^LG_v$
and $\AKvG$ the set of continuous irreducible admissible
representations of $G_v$. We would hope to have a bijection
\begin{equation}
\LKG \lrt \AKG,\ \ \quad \phi\mapsto \pi_\phi.
\end{equation}
Moreover the set ${\mathfrak L}^0_k(G)$ of equivalence classes of
irreducible representations in $\LKG$ should be in bijective
correspondence with the set ${\mathcal A}^0_k(G)$ of all {\it
cuspidal} automorphic representations in $\AKG.$ This would be
supplemented by local bijection
\begin{equation}
\LKvG\lrt {\mathcal A}_{k_v}(G)   \quad {\rm for\ any\ place}\ v\
of\ k.
\end{equation}
The local and global bijections should be {\it compatible} in the
sense that for any $\phi:L_k\lrt \LG,$ there is an automorphic
representation $\pi_\phi=\otimes_v \pi_{\phi,v}$ of $\GA$ with the
correspondence $\phi\mapsto \pi_\phi$ such that for each place $v$
of $k$, the restriction $\phi_v=\phi\circ i_v$ of $\phi$ to
$L_{k_v}$ corresponds to the local component $\pi_{\phi,v}$ of
$\pi_\phi.$ Of course, one expects that all of these
correspondences (4.2) and (4.3) would satisfy properties similar
to those in the local Langlands conjecture, e.g., the preservation
of $L$- and $\varepsilon$-factors with twists, etc.

\vskip 0.2cm The local Langlands groups are elementary. They are
given by
\begin{equation*}
L_{k_v}=
\begin{cases}
W_{k_v} & \quad\quad {\rm if}\ v\ {\rm is\ archimedean},\\
W_{k_v}\times SU(2,\BR) & \quad\quad {\rm if}\ v\ {\rm is\
nonarchimedean},
\end{cases}
\end{equation*}
where $W_{k_v}$ is again the Weil group of $k_v$. Thus the local
Langlands group $L_{k_v}$ is a split extension of $W_{k_v}$ by
compact simply connected Lie group. But the hypothetical Langlands
group will be much larger. It would be an infinite fibre product
of nonsplit extension
\begin{equation}
1\lrt K_c \lrt L_c\lrt W_k\lrt 1
\end{equation}
of the Weil group $W_k$ by a compact, semsimple, simply connected
Lie group $K_c$. However one would have to establish something in
order to show that $L_k$ has all the desired properties.

\vskip 0.2cm Grothendieck's conjectural theory of motives
introduces the so-called {\it motivic Galois} group $\MK$, which
is a reductive proalgebraic group over $\BC$ and comes with a
proalgebraic projection $\MK\lrt \G_k$. A {\it motive} of rank $n$
is to be defined as a proalgebraic representation
\begin{equation*}
{\mathbb M}: \MK\lrt GL(n,\BC).
\end{equation*}
We observe that any continuous representation of $\G_k$ pulls back
to $\MK$ and can be viewed as a motive in the above sense. It is
conjectured that the arithmetic information in any motive $\mathbb
M$ is directly related to analytic information from some
automorphic representations of $\GA$. The conjectural theory of
motives also applies to any completion $k_v$ of $k$. It produces a
proalgebraic group ${\mathcal M}_{k_v}$ over $\BC$ that fits into
a commutative diagram
\begin{equation*}
\begin{array}{ccc}
{\mathcal M}_{k_v} & \longrightarrow & {\mathcal M}_k\\
\Big\downarrow & & \Big\downarrow \\
\Gamma_{k_v} & \hookrightarrow & \Gamma_k
\end{array}
\end{equation*}
\noindent of proalgebraic homomorphisms.

\vskip 0.2cm In 1979, Langlands \cite{Lan7} speculated the
following:

\vskip 0.2cm\noindent {\bf Conjecture B\,(Langlands \cite[Section
2]{Lan7}).} There is a commutative diagram

\begin{equation*}
\begin{array}{ccc}
L_k &\stackrel{\phi}\longrightarrow&  {\mathcal M}_k\\
\ \ \ \,\,\,\,\, \searrow&&\!\!\!\!\!\!\!\swarrow\\
&  \Gamma_k  &\\
\end{array}
\end{equation*}

\noindent together with a compatible commutative diagram

\begin{equation*}
\begin{array}{ccc}
L_{k_v} &\stackrel{\phi_v}\longrightarrow&  {\mathcal M}_{k_v}\\
\ \ \ \,\,\,\,\, \searrow&&\!\!\!\!\!\!\!\swarrow\\
&  \Gamma_k  &\\
\end{array}
\end{equation*}

\noindent for each completion $k_v$ of $k$, in which $\phi$ and
$\phi_v$ are continuous homeomorphisms. There should be the
analogue of the notion of {\it admissibility} of the maps $\phi$
and $\phi_v$ as in the local case (cf.\,\cite[p.\,40]{B}).

\vskip 0.3cm The above conjecture implies that we can attach to
any proalgebraic homomorphism $\mu$ from $\MK$ to $\LG$ over
$\G_k$, its associated automorphic representation of $\GA.$ In
particular, if we take $G=GL(n)$, it means that we can attach to
any motive $\mathbb M$ of rank $n$, an automorphic representation
$\pi_{\mathbb M}=\otimes_v \pi_{\mathbb M,v}$ of $GL(n,\A)$ with
the following property: The family of semi-simple conjugacy
classes $c(\pi_{\mathbb M})=\{ c(\pi_{\mathbb M,v})\}$ in
$GL(n,\BC)$ associated to $\pi_{\mathbb M}$ is equal to the family
of conjugacy classes $c(\mathbb M)=\{ c_v(\mathbb M)\}$ obtained
from $\mathbb M$, and the local homomorphism ${\mathcal
M}_{k_v}\lrt {\mathcal M}_k$ at places $v$ that are unramified for
$\mathbb M$. In fact, $c_v(\mathbb M)$ is the image of the
Frobenius class $Fr_v$ under a different kind of $\G_k$, namely a
compatible family
$$\G_k\lrt \prod_{\ell\notin S(\mathbb M)\cup \{v\}}GL(n,
{\overline \BQ}_{\ell})
$$
of $\ell$-adic representations attached to $\mathbb M$. Our task
now is to find some natural ways to construct an explicit
candidate for ${\mathcal M}_k$ and then to clarify the structure
of ${\mathcal M}_k.$ It is suggested by experts \cite{Ram1} that
${\mathcal M}_k$ be a proalgebraic fibre product of certain
extensions

\begin{equation}
1 \lrt {\mathcal D}_c \lrt {\mathcal M}_c \lrt {\mathcal T}_k \lrt
1
\end{equation}
of a fixed group ${\mathcal T}_k$ by complex, semisimple, simply
connected algebraic groups ${\mathcal D}_c.$ The group ${\mathcal
T}_k$ is an extension
\begin{equation}
1 \lrt {\mathcal S}_k \lrt {\mathcal T}_k \lrt \G_k \lrt 1
\end{equation}
of $\G_k$ by a complex proalgebraic torus ${\mathcal
S}_k.$\,(cf.\,\cite[Chapter II]{Se1},\,\,\cite[Section
5]{Lan7},\,\,\cite[Section 7]{Se2}) The contribution to ${\mathcal
M}_k$ of any ${\mathcal M}_c$ is required to match the
contribution to $L_k$ of a corresponding $L_c$, in which $K_c$ is
a compact real form of ${\mathcal D}_c$. This construction should
have to come with the following diagram

\begin{equation*}
\begin{array}{ccccccccc}
1 & \longrightarrow & L_k^c         & \longrightarrow & L_k & \longrightarrow & \Gamma_k &  \longrightarrow & 1\\
{} & {}             &\Big\downarrow & {}              & \Big\downarrow & {} &\Big\downarrow & {} & {}\\
1 & \longrightarrow & {\mathcal M}_k^c & \longrightarrow &
{\mathcal M} _k & \longrightarrow & \Gamma_k &  \longrightarrow &
1
\end{array}
\end{equation*}
\vskip 0.5cm \noindent where ${\mathcal M}_k^c$ is a complex
pro-reductive group.

\vskip 0.2cm Let $\Psi_k(G)$ be the set of equivalence classes of
continuous, completely reducible homomorphisms of ${\mathcal M}_k$
into $\LG$, and for each place $v$ of $k$, let $\Psi_{k_v}(G)$ be
the set of equivalence classes of continuous, completely reducible
homomorphisms of ${\mathcal M}_{k_v}$ into $\LG_v$. Then one
should have to obtain a bijective correspondence
\begin{equation}
\Psi_k(G) \lrt \AKG.
\end{equation}
This would be supplemented by local bijective correspondences
\begin{equation}
{\mathcal M}_{k_v}(G)\lrt {\mathcal A}_{k_v}(G)
\end{equation}
for all places $v$ of $k$.

\end{section}
\vskip 1cm
%%%%%%%%%%%%%%%%%%%%%%%%%%%%%%%%%%%%%%%%%%%%%%%%%%%%%%%%%%%%%%%%%%%
%
%   Sec 5  Langlands Functoriality
%
%%%%%%%%%%%%%%%%%%%%%%%%%%%%%%%%%%%%%%%%%%%%%%%%%%%%%%%%%%%%%%%%%%%
\begin{section}{{\bf Langlands Functoriality} }
\setcounter{equation}{0}
\newcommand\LH{{}^LH}
\newcommand\GKH{ {\mathcal G}_k(H)}
\newcommand\GKG{ {\mathcal G}_k(G)}

\vskip 0.2cm As we see in section 2, Shahidi \cite{Sha3} gave a
partially affirmative answer to Conjecture A, which is a question
raised by Langlands for a larger class of automorphic
$L$-functions $L(s,\pi,r)$ obtained from cuspidal automorphic
representations $\pi$ of a Levi subgroup $M$ of a quasi-split
reductive group and the adjoint representation of $\LM$ on the
real Lie algebra ${}^L{\mathfrak n}$ of ${}^LN$. This suggest
trying, given a $L$-function and a quasi-split reductive group
$G$, to see whether $G$ has an automorphic representation with the
given $L$-function. Many instances of such questions can be
regarded as special cases of the {\it lifting problem}, nowadays
called the {\it Principle of Functoriality}, with respect to
morphisms of $L$-groups. The motivation of this problem stems from
a global side. There is also a local version for this problem.
These questions were raised and also formulated by Langlands
\cite{Lan1} in the late 1960s. Roughly speaking, the principle of
functoriality describes profound relationships among automorphic
forms on different groups.

\vskip 0.3cm Let $k$ be a local or global field, and let $H,G$ two
connected reductive groups defined over $k$. A homomorphism
$\s:\LH\lrt \LG$ is said to an $L$-{\it homomorphism} if it
satisfies the following conditions (1)-(3):

\vskip 0.2cm (1) $\s$ is a homomorphism over the absolute Galois
group $\G_k$, namely, the following diagram is commutative;
\begin{equation*}
\begin{array}{ccc}
{}^LH  &\stackrel{\sigma}\longrightarrow&  {}^LG\\
\ \ \ \,\,\,\,\, \searrow&&\!\!\!\!\!\!\!\swarrow\\
&  \Gamma_k  &\\
\end{array}
\end{equation*}

(2) $\sigma$ is continuous;\\
\indent (3) The restriction of $\s$ to $\LH^0$ is a complex
analytic homomorphism of $\LH^0$ into $\LG^0.$

\vskip 0.2cm Let $\GKH$\,(resp.\ $\GKG$) be the set of all
admissible homomorphisms $\phi:W_k'\lrt \LH$ (resp.\,$\LG$) modulo
inner automorphisms by elements of $\LH^0$ (resp.\,$\LG^0$).
Suppose $G$ is quasi-split. Given a fixed $L$-homomorphism
$\s:\LH\lrt \LG,$ if $\phi$ is any element in $\GKH$, then the
composition $\s\circ \phi$ is an element in $\GKG.$ It is easily
seen that the correspondence $\phi\mapsto \s\circ \phi$ yields the
canonical map
\begin{equation*}
{\mathcal G}_k(\s):\GKH\lrt \GKG.
\end{equation*}
If $k$ is a global field and $v$ is a place of $k$, then $\LG_v$
can be viewed as a subgroup of $\LG$ because $\G_{k_v}$ is
regarded as a subgroup of $\G_k$. Therefore $\sigma$ induces the
$L$-homomorphism $\s_v:\LH_v\lrt \LG_v$ and hence a local map
\begin{equation*}
{\mathcal G}_k(\s_v):{\mathcal G}_k(H_v)\lrt {\mathcal G}_k(G_v).
\end{equation*}
We refer to \cite[pp. 54-58]{B} for more detail on these stuffs.

\vskip 0.3cm \noindent {\bf Langlands Functoriality Conjecture
(Langlands \cite{Lan1}).} Let $k$ be a global field, and let $H,G$
two connected reductive groups over $k$ with $G$ quasi-split.
Suppose $\s:\LH\lrt \LG$ is an $L$-homomorphism. Then for any
automorphic representation $\pi=\otimes_v \pi_v$ of $H(\A)$, there
exists an automorphic representation $\Pi=\otimes_v \Pi_v$ of
$G(\A)$ such that
\begin{equation}
c(\Pi_v)=\s(c(\pi_v)),\quad v\notin S(\pi)\cup S(\Pi),
\end{equation}
where $S(\pi)$\,(resp. $S(\Pi)$) denotes a finite set of ramified
places of $k$ for $\pi$\,(resp. $\Pi$) so that $\pi_v$\,(resp.
$\Pi_v$) is unramified for every place $v\notin S(\pi)$\,(resp.
$v\notin S(\Pi)).$ We note that the condition (5.1) is equivalent
to the condition
\begin{equation}
L_S(s,\Pi,r)=L_S(s,\pi,r\circ \s),\quad S=S(\pi)\cup S(\Pi)
\end{equation}
for every finite dimensional complex representation $r$ of $\LG$.

\vskip 0.3cm\noindent {\bf Remark 5.1.} For a nonarchimedean local
field $k$, we can formulate a local version of Langlands
Functoriality Conjecture replacing the word ``automorphic" by
``admissible" and modifying some facts of an $L$-homomorphism.

\vskip 0.3cm\noindent {\bf Remark 5.2.} Suppose $k$ is a
nonarchimedean local field with the ring of integers ${\mathcal
O}.$ Suppose $H$ and $G$ are quasi-split and there is a finite
extension $E$ of $k$ such that both $H$ and $G$ split over $E$,
and have an $\mathcal O$-structure so that both $H(\mathcal O)$
and $G(\mathcal O)$ are special maximal compact subgroups. Let
$\pi$ be an unramified representation of $H(k)$ in the sense that
$\pi$ has a nonzero $H(\mathcal O)$-fixed vector, and let
$\phi=\phi_\pi \in \GKH$ be the unramified parameter of $\pi$.
Then for any $L$-homomorphism $\s:\LH\lrt \LG$, the parameter
${\widetilde \phi}=\sigma\circ \phi$ is unramified and defines an
$L$-packet $\prod_{\widetilde \phi}(G)$ which contains exactly one
unramified representation $\Pi$ of $G(k)$ to be called the {\it
natural lift} of $\pi$.(cf. \cite[p. 55]{B})

\vskip 0.3cm If we assume that the Local Langlands Conjecture,
briefly [LLC] is valid, we can reformulate the Langlands
Functoriality Conjecture using [LLC] in the following way. Let
$\pi=\otimes_v \pi_v$ be an automorphic representation of $H(\A)$.
According to [LLC], we can attach to each $\pi_v$, an element
$\phi_v:W_k'\lrt \LH_v$ in ${\mathcal G}_{k_v}(H).$ The
composition $\s\circ \phi_v$ is an element in ${\mathcal
G}_{k_v}(G).$ By [LLC] again, one has an irreducible admissible
representation $\Pi_v$ of $G_v$ attached to $\s\circ \phi_v.$ Then
$\Pi=\otimes_v \Pi_v$ is an irreducible admissible representation
of $G(\A)$. Therefore Langlands Functoriality Conjecture is
equivalent to the statement that $\Pi$ must be an automorphic
representation of $G(\A).$

\vskip 0.2cm If we assume that Global Langlands Conjecture,
briefly [GLC] is valid, we can also reformulate Langlands
Functoriality Conjecture using [GLC] as follows: Given an
automorphic representation $\pi$ of $H(\A)$ with its associated
parameter $\phi_\pi:L_k\lrt \LH,$ there must be an $L$-packet
$\Pi_{\s\circ  \phi_\pi}(G)$ attached to $\s\circ  \phi_\pi.$

\vskip 0.3cm\noindent {\bf Example 5.3.} Suppose $H=\{1\}$ and
$G=GL(n).$ Clearly an automorphic representation $\pi$ of $H(\A)$
is trivial. The choice of $\s$ amounts to that of an admissible
homomorphism
$$ \s: {\rm Gal}(E/k)\lrt GL(n,\BC)=\,\LG$$
for a finite Galois extension $E$ of $k$. Therefore Langlands
Functoriality Conjecture reduces to the following assertion.

\vskip 0.2cm \noindent {\bf Strong Artin Conjecture (Langlands
\cite{Lan1}).} Let $k$ be a number field. For an $n$-dimensional
complex representation $\s$ of ${\rm Gal}(E/k)$, there is an
automorphic representation $\pi$ of $GL(n,\A)$ such that
\begin{equation*}
c(\pi_v)= \s (Fr_v),\quad v\notin S(E),
\end{equation*}
where $S(E)$ is a finite set of places including all the ramified
places of $E$.

\vskip 0.3cm The above conjecture was established partially but
remains unsettled for the most part. We summarize the cases that
have been established until now chronically.

\vskip 0.2cm \noindent (a) The case $n=1$\,: This is the Artin
reciprocity law, namely, $k^* \cong W_k^{\rm ab}$, which is the
essential part of abelian class field theory. The image of $\s$ is
of cyclic type or of dihedral type.\\
(b) The case where $n=2$ and ${\rm Gal}(E/k)$ is solvable\,: The
conjecture was solved by Langlands \cite{Lan8} when the image of
$\s$ in $PSL(2,\BC)$ is of tetrahedral type, that is, isomorphic
to $A_4$, and by Tunnell \cite{Tu1} when the image of $\s$ in
$PSL(2,\BC)$ is of octahedral type, i.e., isomorphic to $S_4$. It
is a consequence of cyclic base change for $GL(2).$ These cases
were used by A. Wiles \cite{Wil} in his proof of Fermat's Last
Theorem.
\\
(c) The case where $n$ is arbitrary and ${\rm Gal}(E/k)$ is
nilpotent\,: The conjecture was established by Arthur and Clozel
\cite{AC} as an application of cyclic base change for $GL(n)$.\\
(d) The case where $n=2$ and the image of $\s$ is of icosahedral
type\,:
Partial results were obtained by Taylor et al.(cf. \cite{BDST}) C. Khare proved this case.\\
(e) The case where $n=4$ and ${\rm Gal}(E/k)$ is solvable\,: The
conjecture was established by Ramakrishnan \cite{Ram3} for
representations $\sigma$ that factor through the group $GO(4,\BC)$
of orthogonal similitudes.

\vskip 0.3cm\noindent {\bf Example 5.4.} Let $k$ be a number
field. Let $H=Sp(2n),\, SO(2n+1),\, SO(2n)$ be the split form, and
$G=GL(N)$, where $N=2n+1$ or $2n$. Then
${}^LSp(2n)=SO(2n+1,\BC),\,{}^LSO(2n+1)=Sp(2n,\BC),\,{}^LSO(2n)=SO(2n,\BC)$
and ${}^LGL(N)=GL(N,\BC).$ As an $L$-homomorphism $\s:\LH\lrt
\LG$, we take the embeddings
\begin{equation*}
{}^LSp(2n)\hookrightarrow GL(2n+1,\BC),\quad {}^LSO(2n+1)
\hookrightarrow GL(2n,\BC),\quad {}^LSO(2n) \hookrightarrow
GL(2n,\BC).
\end{equation*}
In each of these cases, the Langlands weak functorial lift for
irreducible {\it generic} cuspidal automorphic representations of
$H(\A)$ was established by Cogdell, Kim, Piateski-Shapiro and
Shahidi \cite{CKPS1, CKPS2}. Here the notion of ``weak" automorphy
means that an automorphic representation of $GL(n)$ exists whose
automorphic $L$-function matches the desired Euler product except
for a finite number of factors. The proof is based on the converse
theorems for $GL(n)$ established by Cogdell and Piateski-Shapiro
\cite{CP1, CP2}. It is still an open problem to establish the
Langlands functorial lift from irreducible {\it non-generic}
cuspidal automorphic representations of $H(\A)$ to $G(\A)$.

Let $H={\rm GSpin}_m$ be the general spin group of semisimple rank
$[{\frac m2}]$, i.e., a group whose derived group is ${\rm
Spin}_m$. Then the $L$-group of $G$ is given by
$${}^L{\rm GSpin}_m=\begin{cases} GSO_m  \quad & {\rm if}\ m\ {\rm is\ even};\\
GSp_{2[{\frac m2}]} \quad & {\rm if}\ m\ {\rm is\ odd}.
\end{cases}$$
In each case we have an embedding
\begin{equation}
i:{}^LH^0\ \hookrightarrow \ GL(N,\BC),\quad N=m\ {\rm or}\
2\Big[{\frac m2}\Big].
\end{equation}
Asgari and Shahidi \cite{AgSh1, AgSh2} proved that if $\pi$ is a
generic cuspidal representation of ${\rm GSpin}_m$, then the
functoriality is valid for the embedding (5.3).

If $H=SO(2n+1),$ for generic cuspidal representations, Jiang and
Soudry \cite{JiS1} proved that the Langlands functorial lift from
$SO(2n+1)$ to $GL(2n)$ is {\it injective} up to isomorphism. Using
the functorial lifting from $SO(2n+1)$ to $GL(2n)$, Khare, Larsen
and Savin \cite{KLS1} proved that for any prime $\ell$ and any
even positive integer $n$, there are infinitely many exponents $k$
for which the finite simple group $PSp_n({\mathbb F}_{\ell^k})$
appears as a Galois group over $\BQ$. Furthermore, in their recent
paper \cite{KLS2} they extended their earlier work to prove that
for a positive integer $t$, assuming that $t$ is even if $\ell=3$
in the first case (1) below, the following statements (1)-(3)
hold: \vskip 0.15cm (1) Let $\ell$ be a prime. Then there exists
an integer $k$ divisible by $t$ such that the simple group
$G_2({\mathbb F}_{\ell^k})$ appears as a Galois group over
$\BQ$.\\
\indent (2) Let $\ell$ be an odd prime. Then there exists an
integer $k$ divisible by $t$ such that the simple finite group
$SO_{2n+1}({\mathbb F}_{\ell^k})^{\rm der}$  or the finite
classical group $SO_{2n+1}({\mathbb F}_{\ell^k})$ appears as a
Galois group over
$\BQ$.\\
\indent (3) If $\ell\equiv 3,5\, ({\rm mod}\,8)$ and $\ell$ is a
prime, then there exists an integer $k$ divisible by $t$ such that
the simple finite group $SO_{2n+1}({\mathbb F}_{\ell^k})^{\rm
der}$ appears as a Galois group over $\BQ$.

\vskip 0.2cm The construction of Galois groups in (1)-(3) is based
on the functorial lift from $Sp(2n)$ to $GL(2n+1)$, and the
backward lift from $GL(2n+1)$ to $Sp(2n)$ plus the theta lift from
$G_2$ to $Sp(6)$.

\vskip 0.3cm\noindent {\bf Example 5.5.} Let $k$ be a number
field. For two positive integers $m$ and $n$, we let
\begin{equation*}
H=GL(m)\times GL(n)\qquad  {\rm and}\qquad G=GL(mn).
\end{equation*}
Then $\LH=GL(m,\BC)\times GL(n,\BC)$ and $\LG=GL(mn,\BC).$ We take
the $L$-homomorphism
\begin{equation*}
\s:GL(m)\times GL(n)\lrt GL(mn,\BC)
\end{equation*}
given by the tensor product. Suppose $\pi=\otimes_{v} \pi_v$ and
$\tau=\otimes_v \tau_v$ are two cuspidal automorphic
representations of $GL(m,\A)$ and $GL(n,\A)$ respectively. By
[LLC] for $GL(N)$ \cite{HT,Hen2, Lan10}, one has the
parametrizations
\begin{equation*}\phi_v:W_k'\lrt GL(m,\BC)\qquad {\rm and}\qquad
\psi_v:W_k'\lrt GL(n,\BC).
\end{equation*}
Let
$$[\phi_v,\psi_v]:W_k' \lrt GL(m,\BC)\times GL(n,\BC)\hookrightarrow GL(mn,\BC)=\LG$$
be the admissible homomorphism of $W_k'$ into $\LH$ defined by
$$[\phi_v,\psi_v](w)=\big(\phi_v(w),\psi_v(w)\big),\quad w\in W_k'.$$
The composition $\theta_v= \s\circ [\phi_v,\psi_v]$ is an
admissible homomorphism of $W_k'$ into $\LG$. According to [LLC]
for $GL(N)$, one has an irreducible admissible representation of
$GL(mn,k_v)$ attached to $\theta_v$, denoted by $\pi_v \boxtimes
\psi_v$. We set
$$\pi\boxtimes \tau= \bigotimes_v \pi_v \boxtimes
\psi_v.$$ The validity of the Langlands Functoriality Conjecture
for the $L$-homomorphism $\s:\LH\lrt\LG$ implies that
$\pi\boxtimes \tau$ is an automorphic representation of
$GL(mn,\A).$ Ramakrishnan \cite{Ram2} used the converse theorem
for $GL(4)$ of Cogdell and Piatetski-Shapiro to establish the
functoriality for $GL(2)\times GL(2).$ Kim and Shahidi \cite{KiS2}
established the functoriality for $GL(2)\times GL(3)$.

\vskip 0.3cm\noindent {\bf Example 5.6.} Let $H=GL(2)$. For a
positive integer $m\geq 2,$ let $G=GL(m+1).$ Suppose
$\pi=\otimes_v \pi_v$ is an automorphic representation of $H(\A).$
According to [LLC] for $GL(n)$ \cite{HT,Hen2, Lan10}, for each
place $v$ of $k$, we have a semisimple conjugacy class
$c(\pi_v)=\{ {\rm diag}(\alpha_v,\beta_v) \}\subset GL(2,\BC).$ By
[LLC] for $GL(n)$ again, for each place $v$ of $k$, there is an
irreducible admissible representation of $GL(m+1,k_v)$, denoted by
${\rm Sym}^m(\pi_v)$ attached to the semisimple conjugacy class
$$\big\{ {\rm
diag}\big(\alpha_v^m,\alpha_v^{m-1}\beta_v,\cdots,\beta_v^m\big)
\big\}\subset GL(m+1,\BC).$$ We set
\begin{equation*}
{\rm Sym}^m(\pi):=\bigotimes_v {\rm Sym}^m(\pi_v).
\end{equation*}
Then ${\rm Sym}^m(\pi)$ is an irreducible admissible
representation of $GL(m+1,\A).$ The validity of the Langlands
Functoriality Conjecture for the $L$-homomorphism ${\rm
Sym}^m:GL(2,\BC)\lrt GL(m+1,\BC)$ implies that ${\rm Sym}^m(\pi)$
is an automorphic representation of $GL(m+1,\A).$ As a
consequence, we obtain a complete resolution of the
Ramanujan-Petersson conjecture for Maass forms, the Selberg
conjecture for eigenvalues and the Sato-Tate conjecture. In 1978,
Gelbart and Jacquet \cite{GeJ} established the functoriality for
${\rm Sym}^2$ using the converse theorem on $GL(3)$. In 2002, Kim
and Shahidi \cite{KiS2} established the functoriality for ${\rm
Sym}^3$ using the functoriality for $GL(2)\times GL(3)$.
Thereafter Kim \cite{Ki1} established the functoriality for ${\rm
Sym}^4$. The proof is based on the converse theorems for $GL(n)$
established by Cogdell and Piateski-Shapiro \cite{CP1, CP2}. We
refer to \cite{KiS3} for more results on this topic.

\vskip 0.3cm\noindent {\bf Example 5.7.} For a positive integer
$n\geq 2$, we let
\begin{equation*}
H=GL(n) \qquad {\rm and}\qquad G=GL(N),\quad N={{(n-1)n}\over 2}.
\end{equation*}
Let
\begin{equation*}
\wedge^2: \LH=GL(n,\BC)\lrt \LG=GL(N,\BC)
\end{equation*}
be the $L$-homomorphism of $\LH$ into $\LG$ given by the exterior
square. Suppose $\pi=\otimes_v\pi_v$ is a cuspidal automorphic
representation of $GL(n,\A).$ According to [LLC] for $GL(m)$, for
each place $v$ of $k$, one has the admissible homomorphism
$$\phi_v: W_k'\lrt \LH=GL(n,\BC)$$
parameterizing $\pi_v$. The composition $\psi_v=\wedge^2\circ
\phi_v$ again yields an irreducible admissible representation
$\wedge^2 \pi_v$ of $GL(N,k_v)$ for every unramified
representation $\pi_v$. We set
$$\wedge^2\pi=\bigotimes_v \wedge^2 \pi_v.
$$
Then $\wedge^2\pi$ is an irreducible admissible representation of
$GL(N,\A)$. The validity of the Langlands functoriality implies
that $\wedge^2\pi$ is an automorphic representation of $GL(N,\A).$
Kim \cite{Ki1} established the functoriality for the case $n=4$,
that is a functorial lift from $GL(4)$ to $GL(6).$

\vskip 0.3cm\noindent {\bf Remark 5.8.} As we see in Example 5.3,
5.4 and 5.5, the converse theorem for $GL(n)$ obtained by Cogdell
and Piateski-Shapiro plays a crucial role in establishing the
functoriality for $GL(2)\times GL(3),\ {\rm Sym}^3$ and ${\rm
Sym}^4$. There are widely known three methods in establishing the
Langlands functoriality which are based on the theory of the
Selberg-Arthur trace formula \cite{Art2, AC, Lan8, Lan9}, the
converse theorems for $GL(n)$ \cite{CP1, CP2, GeJ} and the theta
correspondence or theta lifting method (R. Howe, J.-S. Li, S.
Kudla et al.).

\vskip 0.2cm According to the above examples and the converse
theorems for $GL(n)$, we see that the importance of the Langlands
Functoriality Conjecture is that automorphic $L$-functions of any
automorphic representations of any group should be the
$L$-functions of automorphic representations of $GL(n,\A).$ In
this sense we can say that $GL(n,\A)$ is speculated to be the
mother of all automorphic representations, and their offspring
$L$-functions are already to have meromorphic continuations and
the standard functional equation.

%\vskip 0.3cm\noindent {\bf Example 5.9.}

%\vskip 0.3cm\noindent {\bf Example 5.10.}

\end{section}

\vskip 1cm
%%%%%%%%%%%%%%%%%%%%%%%%%%%%%%%%%%%%%%%%%%%%%%%%%%%%%%%%%%%%%%%%%%%
%
%   Sec 5  Appendix\:
%
%%%%%%%%%%%%%%%%%%%%%%%%%%%%%%%%%%%%%%%%%%%%%%%%%%%%%%%%%%%%%%%%%%%
\begin{center}{{\large\bf Appendix\,: Converse Theorems}}
\setcounter{equation}{0}
\end{center}

\vskip 0.3cm We have seen that the converse theorems have been
effectively applied to the establishment of the Langlands
functoriality in certain special interesting cases\,(cf. Example
5.4,\,5.5,\,5.6 and 5.7). We understand that the converse theorems
give a criterion for a given irreducible representation of
$GL(n,\A)$ to be {\it automorphic} in terms of the analytic
properties of its associated automorphic $L$-functions. In this
appendix, we give a brief survey of the history of the converse
theorems and survey the recent results in the local converse
theorems.

\vskip 0.2cm The first converse theorem was established by
Hamburger \cite{Ham} in 1921. This theorem states that any
Dirichlet series satisfying the functional equation of the Riemann
zeta function $\zeta(s)$ and suitable regularity conditions must
be a multiple of $\zeta(s)$. More precisely, this theorem can be
formulated:

\vskip 0.2cm\noindent {\bf Theorem A (Hamburger \cite{Ham},
1921).} Let two Dirichlet series $g(s)=\sum_{n\geq 1}a_n n^{-s}$
and $h(s)=\sum_{n\geq 1}b_n n^{-s}$ converge absolutely for ${\rm
Re}(s)>1.$ Suppose that both $(s-1)h(s)$ and $(s-1)g(s)$ are
entire functions of finite order. Assume we have the following
functional equation:
\begin{equation*}
\pi^{-{s\over 2}}\,\Gamma(s/2)\,h(s)=\pi^{-{{1-s}\over
2}}\,\Gamma((1-s)/2)\, g(1-s).
\end{equation*}
Then $g(s)=h(s)=a_1 \zeta(s).$ Here $\Gamma(s)$ is the usual Gamma
function, and an entire function $f(s)$ is said to be {\it of
order} $\rho$ if
$$f(s)=O\big(|s|^{\rho+\epsilon}\big)\qquad {\rm for\ any}\ \epsilon >
0.$$

\vskip 0.3cm Unfortunately Hamburger's converse theorem was not
well recognized until the generalization to $L$-functions attached
to holomorphic modular forms was done by Hecke \cite{Hec} in 1936.
Hecke proved his converse theorem connecting certain $L$-functions
which satisfy a certain functional equation with holomorphic
modular forms with respect to the full modular group $SL(2,\BZ)$.
For a good understanding of Hecke's converse theorem, we need to
describe Hecke's idea and argument roughly. Let
$$f(\tau)=\sum_{n\geq 1} a_n e^{2\pi i n\tau}$$
be a holomorphic modular form of weight $d$ with respect to
$SL(2,\BZ).$ To such a function $f$ Hecke attached an $L$-function
$L(s,f)$ via the Mellin transform
\begin{equation*}
(2\pi)^{-s}\,\Gamma(s)\,L(s,f)=\int_0^{\infty}
f(iy)\,y^s\,{{dy}\over y}
\end{equation*}
and derived the functional equation for $L(s,f).$ He inverted this
process by taking a Dirichlet series
$$D(s)=\sum_{n\geq 1} {{a_n}\over {n^s}}$$
and assuming that it converges absolutely in some half plane, has
an analytic continuation to an entire function of finite order,
and satisfies the same functional equation as $L(s,f).$ In his
masterpiece \cite{Hec}, he could reconstruct a holomorphic modular
form from $D(s)$ by Mellin inversion
\begin{equation*}
f(iy)=\sum_{n\geq 1}a_ne^{-2\pi ny}={1\over{2\pi i}}
\int_{2-i\infty}^{2+i\infty}(2\pi)^{-s}\,\Gamma(s)\,D(y)\,y^s\,ds
\end{equation*}
and obtain the modular transformation law for $f(\tau)$ under
$\tau\mapsto -\tau^{-1}$ from the functional equation for $D(s)$
under $s\mapsto d-s$. This is Hecke's converse theorem\,! You
might agree that Hecke's original idea and argument are remarkably
beautiful. In 1949, in his seminal paper \cite{Ma}, Maass, a
student of Hecke, extended Hecke's method to non-holomorphic forms
for $SL(2,\BZ).$ In 1967, the next very important step was made by
Weil in his paper \cite{We} dedicated to C. L. Siegel (1896-1981)
celebrating Siegel's seventieth birthday. Weil showed how to work
with Dirichlet series attached to holomorphic modular forms with
respect to congruence subgroups of $SL(2,\BZ).$ He proved that if
a Dirichlet series together with a sufficient number of twists
satisfies {\it nice} analytic properties and functional equations
with reasonably suitable regularity, then it stems from a
holomorphic modular form with respect to a congruence subgroup of
$SL(2,\BZ).$ More precisely his converse theorem can be formulated
as follows.

\vskip 0.2cm\noindent {\bf Theorem B (Weil \cite{We}, 1967).} Fix
two positive integers $d$ and $N$. Suppose the Dirichlet series
$$D(s)=\sum_{n\geq 1} {{a_n}\over {n^s}}$$
satisfies the following properties: \vskip 0.2cm\noindent (W1)
$D(s)$
converges absolutely for sufficiently large ${\rm Re}(s)\gg 0$\,;\\
(W2) For every primitive character $\chi$ of modulus $r$ with
$(r,N)=1$, the function
\begin{equation*}
\Lambda(s,\chi):=(2\pi)^{-s}\,\Gamma(s)\,\sum_{n\geq
1}{{a_n\chi(n)}\over {n^s}}
\end{equation*}
has an analytic continuation to an entire function of $s$ to the
whole complex plane, and is bounded
in vertical strips of finite width;\\
(W3) Every such a function $\Lambda(s,\chi)$ satisfies the
functional equation
\begin{equation*}
\Lambda(s,\chi)=w_\chi\, r^{-1}\,(r^2N)^{{\frac d2}-s}\,
\Lambda(d-s,{\overline\chi}),
\end{equation*}
where
$$w_\chi =i^d\,\chi(N)\,g(\chi)^2$$
and
$$g(\chi)=\sum_{n\,({\rm mod}\ r)}\chi(n)\,e^{2\pi in/r}.$$
Then the function
\begin{equation*}
f(\tau)=\sum_{n\geq 1}a_n e^{2\pi in\tau},\qquad \tau\in {\mathbb
H}=\{ \tau\in\BC\,|\ {\rm Im}\,\tau >0\,\}
\end{equation*}
is a holomorphic cusp form of weight $d$ with respect to the
congruence subgroup $\G_0(N)$, where
$$\Gamma_0(N)=\left\{ \begin{pmatrix} a & b \\ c & d \end{pmatrix}
\in SL(2,\BZ) \Big|\ c\equiv 0\ ({\rm mod}\ N)\  \right\}.$$
Weil's proof follows closely Hecke's idea and argument. We see
that his converse theorem provides a condition for modularity of
the Dirichlet series $D(s)$ under $\G_0(N)$ in terms of the
functional equations of Dirichlet series {\it twisted by primitive
characters}. In fact, Weil's converse theorem influenced the
complete proof of the Shimura-Taniyama conjecture given by Wiles
\cite{Wil}, Taylor et al. So we can say that the work of Weil
marks the beginning of the modern era in the study of the
connection between $L$-functions and automorphic forms.

\vskip 0.2cm In 1970 Jacquet and Langlands \cite{JL} established
the converse theorem for $GL(2)$ in the adelic context of
automorphic representations of $GL(2,\A)$ based on Hecke's
original idea. In 1979 Jacquet, Piatetski-Shapiro and Shalika
\cite{JPS1} established the converse theorem for $GL(3)$ in the
adelic context. Finally in 1994, generalizing the work on the
converse theorems on $GL(2)$ and $GL(3)$, Cogdell and
Piatetski-Shapiro \cite{Co, CP1, CP2} proved the converse theorem
for $GL(n)$ with arbitrary $n\geq 1$ in the context of automorphic
representations. The idea and technique in the proof of Cogdell
and Piatetski-Shapiro are surprisingly almost the same as Hecke's.
We now describe the converse theorems formulated and proved by
them.

\vskip 0.2cm Let $k$ be a global field, $\A$ its adele ring, and
let $\psi$ be a fixed nontrivial continuous additive character of
$\A$ which is trivial on $k$. Let $\pi=\otimes_v\pi_v$ be an
irreducible admissible representation of $GL(n,\A)$, and let
$\tau=\otimes_v\tau_v$ be a cuspidal automorphic representation of
$GL(m,\A)$ with $m<n.$ We define formally
\begin{equation*}
L(s,\pi\times\tau)=\prod_v L(s,\pi_v\times\tau_v)\quad {\rm
and}\quad \varepsilon (s,\pi\times \tau)=\prod_v \varepsilon
(s,\pi_v\times \tau_v,\psi_v).
\end{equation*}
\indent We say that $L(s,\pi\times \tau)$ is {\it nice} if it
satisfies the following properties\,:

\vskip 0.2cm \noindent (N1) $L(s,\pi\times \tau)$ and
$L(s,{\widetilde\pi}\times {\widetilde\tau})$ have analytic
continuations to entire functions, where ${\widetilde\pi}$ (resp. ${\widetilde\tau}$)
denotes the contragredient of $\pi$ (resp.\,$\tau$);\\
(N2) $L(s,\pi\times \tau)$ and $L(s,{\widetilde\pi}\times
{\widetilde\tau})$ are bounded in vertical strips of finite
width;\\
(N3) These entire functions satisfy the standard functional
equation
\begin{equation*}
L(s,\pi\times\tau)=\varepsilon (s,\pi\times
\tau)\,L(1-s,{\widetilde\pi}\times {\widetilde\tau}).
\end{equation*}

\vskip 0.3cm\noindent {\bf Theorem C\,(Cogdell and Piatetski
\cite{CP1, CP2}, 1994).} Let $\pi$ be an irreducible admissible
representation of $GL(n,\A)$ whose central character is trivial on
$k^*$ and whose $L$-function $L(s,\pi)$ converges absolutely in
some half plane. Assume that $L(s,\pi\times \tau)$ is {\it nice}
for every cuspidal automorphic representation $\tau$ of $GL(m,\A)$
for $1\leq m\leq n-2$. Then $\pi$ is a {\it cuspidal} automorphic
repesentation of $GL(n,\A).$

\vskip 0.2cm Furthermore they proved the following theorem.

\vskip 0.3cm\noindent {\bf Theorem D\,(Cogdell and Piatetski
\cite{CP2}, 1999).} Let $\pi$ be an irreducible admissible
representation of $GL(n,\A)$ whose central character is trivial on
$k^*$ and whose $L$-function $L(s,\pi)$ converges absolutely in
some half plane. Let $S$ be a finite set of finite places. Assume
that $L(s,\pi\times \tau)$ is {\it nice} for every cuspidal
automorphic representation $\tau$ of $GL(m,\A)$ for $1\leq m\leq
n-2$, which is unramified at the places in $S$. Then $\pi$ is
quasi-automorphic in the sense that there is an automorphic
representation $\pi'$ of $GL(n,\A)$ such that $\pi_v\cong \pi'_v$
for all $v\notin S.$

\vskip 0.2cm The local converse theorem for $GL(n)$ was first
formulated by Piatetski-Shapiro in his unpublished Maryland notes
(1976) with his idea of deducing the local converse theorem from
his global converse theorem. It was proved by Henniart \cite{Hen1}
using a local approach. The local converse theorem is a basic
ingredient in the proof of [LCC] for $GL(n)$ by Harris and Taylor
\cite{HT} and by Henniart \cite{Hen2}.

\vskip 0.2cm The local converse theorem for $GL(n)$ can be
formulated as follows.

\vskip 0.2cm \noindent {\bf Theorem E (Henniart \cite{Hen1},
1993).} Let $k$ be a nonarchimedean local field of characteristic
$0$. Let $\tau$ and $\tau'$ be irreducible admissible generic
representations of $GL(n,k)$ with the same central character.
Assume the twisted local gamma factors (cf.\,\cite{JPS2}) are the
same, i.e.,
\begin{equation*}
\gamma (s,\tau\times \rho,\psi)=\gamma (s,\tau'\times \rho,\psi)
\end{equation*}
for all irreducible supercuspidal representations $\rho$ of
$GL(m,k)$ with $1\leq m\leq n-1.$ Then $\tau$ is isomorphic to
$\tau'.$

\vskip 0.2cm\noindent {\bf Remark 1.} It is known that the
twisting condition on $m$ reduces from $n-1$ to $n-2$. It is
expected as a conjecture of H. Jacquet\,\cite[Conjecture 8.1]{CP2}
that the twisting condition on $m$ should be reduced from $n-1$ to
$[{\frac n2}].$

\vskip 0.2cm\noindent {\bf Remark 2.} The local converse theorem
for generic representations of $U(2,1)$ and for $GSp(4)$ was
established by E. M. Baruch in his Ph.D. thesis (Yale Univ.,
1995).

\vskip 0.2cm Jiang and Soudry \cite{JiS1} proved the local
converse theorem for irreducible admissible generic
representations of $SO(2n+1,k).$

\vskip 0.2cm\noindent {\bf Theorem F (Jiang and Soudry
\cite{JiS1}, 2003).} Let $\s$ and $\s'$ be irreducible admissible
generic representations of $SO(2n+1,k)$. Assume the twisted local
gamma factors are the same, i.e.
$$\gamma(s,\s\times \rho,\psi)=\gamma(s,\s'\times \rho,\psi)$$
for all irreducible supercuspidal representations $\rho$ of
$GL(m,k)$ with $1\leq m\leq 2n-1.$ Then $\s$ is isomorphic to
$\s'.$

\vskip 0.2cm I shall give a brief sketch of the idea of their
proof. They first reduce the proof of Theorem F to the case where
both $\s$ and $\s'$ are supercuspidal by studying the existence of
poles of twisted local gamma factors and related properties.
Developing the explicit local Howe duality for irreducible
admissible generic representations of $SO(2n+1,k)$ and the
metaplectic group ${\widetilde{Sp}}(2n,k),$ and using the global
weak Langlands functorial lifting form $SO(2n+1)$ to
$GL(2n)$\,(cf.\,Example 5.4,\,\cite{CKPS1, CKPS2}) and the local
backward lifting from $GL(2n,k)$ to ${\widetilde{Sp}}(2n,k),$ they
relate the local converse theorem for $SO(2n+1)$ with that for
$GL(2n)$ which is well known now.

\vskip 0.2cm As an application of Theorem F, I repeat again that
Jiang and Soudry \cite{JiS1, JiS2} proved the Local Langlands
Reciprocity Law for $SO(2n+1)$. More precisely, there exists a
{\it unique} bijective correspondence between the set of conjugacy
classes of all $2n$-dimensional, admissible, completely reducible,
multiplicity-free, symplectic complex representations of the Weil
group $W_k$ and the set of all equivalence classes of irreducible
generic supercuspidal representations of $SO(2n+1,k)$, which
preserves the relevant local factors. As an application of Theorem
$F$ to the global theory, they proved that the weak Langlands
functorial lifting from irreducible {\it generic} cuspidal
automorphic representations of $SO(2n+1)$ to irreducible
automorphic representations of $GL(2n)$ is {\it injective} up to
isomorphism. It is still an open problem to establish the
Langlands functorial lift from irreducible {\it non-generic}
cuspidal automorphic representations of $SO(2n+1)$ to $GL(2n).$ As
another application to the global theory, they proved the {\it
rigidity theorem} in the sense that if $\pi=\otimes_v\pi_v$ and
$\tau=\otimes_v\tau_v$ are irreducible generic cuspidal
automorphic representations of $SO(2n+1,\A)$ such that $\pi_v$ is
isomorphic to $\tau_v$ for almost all local places $v$, then $\pi$
is isomorphic to $\tau.$

%%%%%%%%%%%%%%%%%%%%%%%%%%%%%%%%%%%%%%%%%%%%%%%%%%%%%%%%%%%%%%%%%%%
%
%   Sec 5  Siegel Modular Forms
%
%%%%%%%%%%%%%%%%%%%%%%%%%%%%%%%%%%%%%%%%%%%%%%%%%%%%%%%%%%%%%%%%%%%
%\begin{section}{{\bf Siegel Modular Forms}
%\setcounter{equation}{0}
%\end{section}

\end{document}